\documentclass[10pt]{amsart}
\usepackage{amsfonts}
\usepackage{amsmath}
\usepackage{amsthm}
\usepackage{amssymb}
\usepackage{latexsym}
\usepackage{multicol}
\usepackage{verbatim}

\advance\textwidth by 1.2in \advance\oddsidemargin by -.6in
\advance\evensidemargin by -.6in
\newtheorem*{cor}{Corollary}
\newtheorem*{lem}{Lemma}
\newtheorem*{prop}{Proposition}

\theoremstyle{definition}
\newtheorem*{defn}{Definition}
\theoremstyle{definition}
\newtheorem*{thm}{Theorem}
\newtheorem*{conj}{Conjecture}
\newtheorem*{rem}{Remark}

\newenvironment{pf}{\proof}{\endproof}
\newcounter{cnt}
\newenvironment{enumerit}{\begin{list}{{\hfill\rm(\roman{cnt})\hfill}}{%
\settowidth{\labelwidth}{{\rm(iv)}}\leftmargin=\labelwidth%
\advance\leftmargin by
\labelsep\rightmargin=0pt\usecounter{cnt}}}{\end{list}}

\theoremstyle{remark}

\numberwithin{equation}{section} 

\def\wt{{\rm wt}}
\def\bxj {\underline{\bx}_j}
\def\bxu {\underline{\bx}_r}

\def\opl_#1{\text{\scriptsize$\bigoplus\limits_{\text{\footnotesize$#1$}}$}}

\newcommand{\bin}[2]{\left[{#1}\atop{#2}\right]}

\begin{document}

\newcommand{\thmref}[1]{Theorem~\ref{#1}}
\newcommand{\secref}[1]{Section~\ref{#1}}
\newcommand{\lemref}[1]{Lemma~\ref{#1}}
\newcommand{\propref}[1]{Proposition~\ref{#1}}
\newcommand{\corref}[1]{Corollary~\ref{#1}}
\newcommand{\remref}[1]{Remark~\ref{#1}}
\newcommand{\defref}[1]{Definition~\ref{#1}}
\newcommand{\er}[1]{(\ref{#1})}
\newcommand{\id}{\operatorname{id}}
\newcommand{\tensor}{\otimes}
\newcommand{\nc}{\newcommand}
\newcommand{\rnc}{\renewcommand}
\newcommand{\qbinom}[2]{\genfrac[]{0pt}0{#1}{#2}}
\nc{\cal}{\mathcal} \nc{\goth}{\mathfrak} \rnc{\bold}{\mathbf}
\renewcommand{\frak}{\mathfrak}
\newcommand{\desc}{\operatorname{desc}}
\newcommand{\Maj}{\operatorname{Maj}}
\renewcommand{\Bbb}{\mathbb}
\nc\bxi{{\mbox{\boldmath $\xi$}}} \nc\bvpi{{\mbox{\boldmath
$\varpi$}}}
 \nc\balpha{{\mbox{\boldmath $\alpha$}}}
 \nc\bpi{{\mbox{\boldmath
$\pi$}}}
 \nc\blambda{{\mbox{\boldmath $\lambda$}}}
\newcommand{\lie}[1]{\mathfrak{#1}}
\makeatletter
\def\section{\def\@secnumfont{\mdseries}\@startsection{section}{1}%
  \z@{.7\linespacing\@plus\linespacing}{.5\linespacing}%
  {\normalfont\scshape\centering}}
\def\subsection{\def\@secnumfont{\bfseries}\@startsection{subsection}{2}%
  \z@{.5\linespacing\@plus.7\linespacing}{-.5em}%
  {\normalfont\bfseries}}
\makeatother
\def\subl#1{\subsubsection{}\label{#1}}
\newcommand{\bem}[1]{\left(\begin{array}{#1}}
\newcommand{\enm}{\end{array}\right)}

\nc{\Cal}{\cal} \nc{\Xp}[1]{X^+(#1)} \nc{\Xm}[1]{X^-(#1)}
\nc{\on}{\operatorname} \nc{\ch}{{\rm ch}} \nc{\Z}{{\bold Z}}
\nc{\J}{{\cal J}} \nc{\C}{{\bold C}} \nc{\Q}{{\bold Q}}
\renewcommand{\P}{{\cal P}}
\nc{\N}{{\Bbb N}} \nc\boa{\bold a} \nc\bob{\bold b} \nc\boc{\bold
c} \nc\bod{\bold d} \nc\boe{{\bold e}} \nc\bof{\bold f}
\nc\bog{\bold g} \nc\boh{\bold h} \nc\boi{\bold i} \nc\boj{\bold
j} \nc\bok{\bold k} \nc\bol{\bold l} \nc\bom{\bold m}
\nc\bon{\bold n} \nc\boo{\bold o} \nc\bop{\bold p} \nc\boq{\bold
q} \nc\bor{\bold r} \nc\bos{{\bold s}} \nc\bou{\bold u}
\nc\bov{\bold v} \nc\bow{\bold w} \nc\boz{\bold z}
\nc\bpr{\bop\bor^-} \nc\ba{\bold A} \nc\bb{\bold B} \nc\bc{\bold
C} \nc\bd{\bold D} \nc\be{\bold E} \nc\bg{\bold G} \nc\bh{\bold H}
\nc\bi{\bold I} \nc\bj{\bold J} \nc\bk{\bold K} \nc\bl{\bold L}
\nc\bm{\bold M} \nc\bn{\bold N} \nc\bo{\bold O} \nc\bp{\bold P}
\nc\bq{\bold Q} \nc\br{\bold R} \nc\bs{\bold S} \nc\bt{\bold T}
\nc\bu{\bold U} \nc\bv{\bold V} \nc\bw{\bold W} \nc\bz{\bold Z}
\nc\bx{\bold x} \nc\ul{{\underline{\ell}}}
\nc\us{{\underline{\bos}}}
\nc\ud{{\underline{d}}} \nc\bF{{\bf F}} \nc\gr{{\rm gr}}

\def \eql {\,{\vartriangleright}\,}
\def \eqd {\,{\blacktriangleright}\,}

\def \alt {{\Huge\bf  [*}}
\def \ea  {{\Huge\bf  *]}}
\nc{\dbgn}[1]{\footnote{#1}}

\title{Weyl, Demazure and  Fusion modules \\for the current algebra  of
$\lie{sl}_{r+1}$}
\author{Vyjayanthi Chari and Sergei Loktev}
\address{VC: Department of Mathematics, University of
California, Riverside, CA 92521} \email{chari@math.ucr.edu}
\address{SL:  Institute for Theoretical and
Experimental Physics, Moscow 117218.}
 \email{ loktev@itep.ru}

\begin{abstract}
We construct a Poincare-Birkhoff--Witt type basis for the  Weyl
modules \cite{CPweyl}  of the current algebra of $\lie{sl}_{r+1}$. As a
corollary we prove the
 conjecture made in \cite{CPweyl}, \cite{CP2}
 on the dimension of the Weyl modules in
this case.
 Further, we relate
the Weyl modules to the fusion modules defined in \cite{FL} of the
current algebra and the Demazure modules in level one
representations of the corresponding affine algebra. In
particular, this allows us to establish substantial cases of the
conjectures in \cite{FL}  on the structure and graded character of
the fusion modules.
\end{abstract}

\maketitle
\setcounter{section}{0}

\section*{Introduction} The study of finite--dimensional
representations of quantum affine algebras has been of some
interest in recent years and there are several different
approaches to the subject, \cite{CPweyl},
\cite{FR},\cite{FM},\cite{Nak},\cite{Ka1} for instance.

An approach that was developed in \cite{CPweyl} was to study these
representations by specializing to the case of classical affine
algebras. It was noticed in that paper, that an irreducible
finite--dimensional representation of a quantum affine algebra, in
general specialized to an indecomposable but reducible
representation of the affine Lie algebra. This behavior is  seen
in the representation theory of simple Lie algebras, when passing
from the characteristic zero situation  to the  case of non--zero
characteristic. The modules in non--zero characteristic which are
obtained from the irreducible finite--dimensional modules of the
simple Lie algebra are called {\em Weyl modules} and are  given by
the same  generators and relations as the modules in
characteristic zero. This analogy motivated the definition of the
Weyl modules in \cite{CPweyl}, in terms of generators and
relations for both the classical and quantum affine Lie  algebras.
It was also conjectured there, that the Weyl modules of the affine
algebra  are the  classical ($q \to 1$)  limit of the standard
modules of the quantum affine algebra defined in \cite{Nak}. In
particular, a proof of the conjecture would give generators and
relations for these representations, or equivalently, would prove
that the standard modules are isomorphic to the quantum Weyl
modules. This latter result was proved recently in \cite{CdM}
using some deep geometric results of Nakajima.

However,  the results of \cite{CPweyl} and  \cite{CP2} showed
 that the conjecture on
the isomorphism between the quantum Weyl modules and the standard
modules would also follow, if a further conjecture on the
dimension of the classical Weyl modules could be proved. This
latter conjecture was established in \cite{CPweyl} for
$\lie{sl}_2$. Moreover, it was shown that it was enough to study
the analogous modules for the current algebra of a simple Lie
algebra, namely the natural parabolic subalgebra of the affine Lie
algebra.

In this paper we show that the conjecture on the dimension of the
classical Weyl modules  is true for $ \lie{sl}_{r+1}$, by constructing
an explicit basis for the Weyl modules over the current algebra of
$\lie{sl}_{r+1}$. We use the basis to give  a graded  fermionic type
character formula for the Weyl modules (see also \cite{HKOTT}).
 We then make connections with several other interesting problems in
the representation theory of affine and current algebras. Thus, we
are able to establish a substantial case of conjectures in
\cite{FL} on the fusion modules for a current Lie algebras. The
definition of these  modules was motivated by studying the space
of conformal blocks  for affine Lie algebras.  The fusion modules
are given by a set of representations of the simple Lie algebra
and a set of complex points, one for each representation. It was
conjectured that  in the case of a simple Lie algebra and where
the representations were irreducible and finite--dimensional the
fusion modules are independent of these points.  This conjecture
is established in this paper  for the fusion of fundamental
representation of $\lie{sl}_{r+1}$ by showing that the fusion module is
isomorphic to a Weyl modules.

Another well--studied family of modules are the Demazure modules.
These modules for the current algebra are obtained by taking the
current algebra module generated by the extremal vectors in
modules of positive level of affine Lie algebras.
 The dimension (and,
actually, character)
 of these modules was computed in \cite{KMOTU}, \cite{M}.  We use these
results
 to  prove that the Weyl modules are isomorphic to the Demazure
modules in the level one representations of the affine Lie algebra
of $\lie{sl}_{r+1}$.  Further, we prove that    the graded character of
the fusion modules can be written in terms of Kostka polynomials
as conjectured in \cite{FL}. Moreover, since the Weyl modules for
the current algebra  can be regarded as a pull--back of Weyl
modules for the affine Lie algebra, it follows that the
$\lie{sl}_{r+1}$--structure of the  Demazure modules can be extended to
 a structure of affine Lie algebra modules. This is related
to a conjecture in \cite{FLitt}.

 For an arbitrary simple Lie algebra, the conjecture
appears to be harder to establish. This is not surprising, since
the representation theory of the corresponding quantum affine
algebra is much more complicated. In particular the Weyl module
associated with a fundamental weight is no longer irreducible as a
module for the simple Lie algebra. But we are convinced that the
conjecture of \cite{CPweyl} is true for these modules and can be
established in a similar way.  Note however, that since the
dimension of fusion modules is known from the definition  the
isomorphism between the Weyl modules, and the  corresponding
fusion modules
 would follow as in the case of $\lie{sl}_{r+1}$ once the
conjecture on dimension of the Weyl modules is proved.

The paper is organized as follows. In Section 1, we define the
Weyl modules and recall some results from \cite{CPweyl}. We also
recall the definition of the fusion modules from \cite{FL} and the
Demazure modules and prove that these modules are quotients of
Weyl modules. We then state the main theorem on the dimension of
the Weyl modules and show that it gives a sufficient condition for
an isomorphism to exist between the Weyl modules,  the  fusion
modules and the Demazure modules. In Section 2 we  start proving
the conjecture by identifying a basis for the Weyl module. The
proof that it is actually a basis involves constructing  a
filtration (which we call a {\em Gelfand--Tsetlin filtration}) for
the Weyl modules, studying the associated graded spaces and using
an induction on $r$. Section 3, is devoted to proving that the
associated graded space is actually isomorphic to a sum of Weyl
modules for $\lie{sl}_r$. This latter result is the motivation for
calling it a {\em Gelfand--Tsetlin} filtration.  We provide an
index of notation at the end of the paper.

{\bf Acknowledgments.} We are grateful to B.~Feigin,
P.~Littelmann, T.~Miwa, E.~Mukhin, and M.~Okado for useful
discussions. SL is partially supported by the Grants  CRDF
RM1-2545-MO-03, PICS 2094, RFBR - 04-01-00702, RF President
Grants MK-3419.2005.1 and N.Sh-8004.2006.2. VC is partially
supported by the NSF grant  DMS-0500751.

{\it{Added in Proof: }}  There has been progress  on the various
conjectures since our paper was posted on the web as
math.QA/0502165. Thus, in math.RT/0509276 the authors prove that
the   Demazure modules, the Weyl modules and the fusion modules
are all isomorphic  for simply--laced Lie algebras extending our
results for $sl_{r+1}$. In particular, this establishes the
conjecture on the dimension of the Weyl modules  in \cite{CPweyl}
for the simply laced algebras. For the nonsimply--laced case, they
give an example to show that the Demazure modules are smaller than
the Weyl modules and formulate a natural  conjecture for the
Demazure modules. They have also have a conjecture for the Weyl
modules which would be an immediate consequence of the conjecture
in \cite{CPweyl}. The methods they use are quite different from
ours and in particular they do not have an analog of the fermionic
character  formula that we give in Section 2 of this paper.
\section{The main theorem and it s applications}
\subsection{Preliminaries}

 Throughout the paper we restrict ourselves to the case of the Lie
algebra
 of $(r+1)\times (r+1)$ trace zero matrices. However, we prefer to use
the more general notation of
 simple Lie algebras since we expect the results to hold in that
 generality.

\subl{} Let  $\bz$ denote  the set of integers, $\bn$ the set of
non--negative integers and $\bn_+$ the set of positive integers.
 Let
 $\lie g = \lie{sl}_{r+1}$ be
the Lie algebra of $(r+1)\times (r+1)$--matrices of trace zero,
$\lie h$ be  the
 Cartan subalgebra of $\lie g$ consisting of diagonal matrices and
$\alpha_i$, $1\le i\le r$, a set of simple roots for $\lie g$ with
respect to $\lie h$. For $1\le i\le j\le r$, let
$\alpha_{i,j}=(\alpha_i+\cdots +\alpha_j)$ and let $x_{i,j}^+$
(resp. $x_{i,j}^-$) be the $(r+1)\times (r+1)$--matrix with $1$ in
the $(i,j+1)^{th}$ (resp. $(j+1,i)^{th}$)--position and $0$
elsewhere. Define subalgebras $ \lie n^\pm$ of $\lie g$ by $$\lie
n^\pm=\bigoplus_{1\le i\le j\le r}\bc x_{i,j}^\pm.$$ For $1\le i
\le r+1$ let $H_i$  be the diagonal matrix with $1$ in the
$i^{th}$ place and zero elsewhere. The elements $h_i =H_i -
H_{i+1}$, $1\le i\le r$ are  a basis of $\lie h$. Let $\{\omega_i
: 1\le i\le r\}$ be the the set  of fundamental weights of $\lie
g$, namely the  basis of $\lie h^*$ which is dual to
$\{\alpha_i:1\le i\le r\}$.  Let $P=\sum_{i=1}^r\bz\omega_i$,
(resp. $Q=\sum_{i=1}^r\bz\alpha_i$) be the  weight lattice (resp.
root lattice) of $\lie g$ and set $P^+=\sum_{i=1}^r\bn\omega_i$,
(resp. $Q^+=\sum_{i=1}^r\bn\alpha_i$). Given a Lie algebra $\lie
a$, we let  $U(\lie a)$ denote the universal enveloping algebra
$\lie a$.

\subl{} Let $\bz[P]$ be the integral group ring of $P$ with basis
$e(\mu)$, $\mu\in P$. If $V$ is any finite--dimensional $\lie
g$--module with
\begin{equation*}
V=\bigoplus_{\mu\in P}V_\mu,\ \ V_\mu=\{v\in V: hv=\mu(h)v,\ \
h\in\lie h\},
\end{equation*}
let $\ch_{\lie g}(V)=\sum_{\mu\in P}\dim(V_\mu)e(\mu)\in \bz[P]$ be
the character of $V$.

Given $\lambda=\sum_{i=1}^r m_i\omega_i\in P^+$, let  $V(\lambda)$
denote the irreducible finite--dimensional $\lie g$--module with
highest weight $\lambda$ and highest weight vector $v_\lambda$.
More precisely, $V(\lambda)=U(\lie g)v_\lambda$ with defining
relations:
\begin{equation}\label{defrel} x_{i,j}^+v_\lambda=0,\ \
(h-\lambda(h))v_\lambda=0,\ \ (x_{i,i}^-)^{m_i+1}v_\lambda=0,
\end{equation}
for all $1\le i\le j\le r$ and  $h\in\lie h$.

\subl{}  Let $\bc[t]$ denote the polynomial ring in an
indeterminate $t$ and for any Lie algebra $\lie a$, set $\lie
a[t]=\lie a\otimes \bc[t]$. Clearly, $\lie a[t]$ is a Lie algebra
with the Lie bracket being given by $[x\otimes f, y\otimes
g]=[x,y]\otimes fg$ for all $x,y\in\lie a$ and $f,g\in\bc[t]$. We
regard $\lie a$ as a subalgebra of $\lie a[t]$ by mapping $x\to
x\otimes 1$ for $x\in\lie a$. We denote the maximal ideal of $\lie
a[t]$ generated by elements of the form $a\otimes t^n$, $a\in\lie
a$, $n>0$ by $\lie at[t]$. The Lie algebra $\lie a[t]$ is a
$\bn$--graded Lie algebra, the grading being given by powers of
$t$ and hence  $\bu(\lie a[t])$ is a $\bn$--graded algebra.

  Given a $\bn$--graded $\lie g[t]$--module $M
=\bigoplus_{s\in\bn} M_s$. Observe  that $M_s$, $s\in\bn$ are
$\lie g$--submodules of $M$. In the case when
\begin{equation}\label{wtd}
M =\bigoplus_{(\mu,s)\in P\times\bn} M_{\mu,s},\ \
M_{\mu,s}=\{v\in M_s: hv=\mu(h)v,\ \ h\in\lie h\},
\end{equation}
recall, that the graded character of $M$ to be the element of
$\bz[P][t]$ is defined  as, $$\ch_t(M)=\sum_{(\mu,s)\in P\times
\bn} \dim(M_{\mu,s})t^se(\mu)=\sum_{s\in\bn}\ch_\lie g (M_s)t^s.$$

\subsection{Weyl modules}
\subl{}

We now recall the definition of the  {\em Weyl modules}.

\begin{defn}
Given $\lambda=\sum_{i=1}^rm_i\omega_i\in P^+$,  the Weyl module
 $W(\lambda)$
is the $\lie g[t]$--module generated by an element $w_\lambda$
with defining relations:
\begin{equation}\label{weylgenrel} \lie
n^+[t]w_\lambda=0,\ \ \lie ht[t]w_\lambda =0,\ \
hw_\lambda=\lambda(h)w_\lambda,\ \ (x^-_{i,i}\otimes
1)^{m_i+1}w_\lambda=0,
\end{equation}
for all $h\in\lie h$, $1\le i \le r$.\hfill\qedsymbol
\end{defn}
\begin{rem} The modules $W(\lambda)$ were defined and
initially studied in \cite{CPweyl} for the affine Lie algebras.
All the results of that paper go over with no difficulty to the
case of the current algebras, in fact an inspection of the proofs
given there show that the results for the affine algebras were
obtained by considering the
 action of the current algebras $\lie g \otimes \bc[t]$. (See
Subsection~\ref{subs_w}
 for further
 details.)
The notation used in this paper is different from that in
\cite{CPweyl}, \cite{CP2}. In those papers, the modules for the
affine Lie algebra were denoted as  $W(\bpi)$ where $\bpi$ is an
$r$--tuple of polynomials in an indeterminate $u$ with constant
term one. For $\lambda=\sum_{i=1}^r m_i\omega_i\in P^+$ and
$a\in\bc^\times$, let $$\bpi_\lambda = ((1-au)^{m_1},\cdots
,(1-au)^{m_r}).$$ The module $W(\lambda)$ is obtained from
$W(\bpi_\lambda)$ by pulling back $W(\bpi_\lambda)$ by the Lie
algebra homomorphism $\lie g[t]\to\lie g[t]$ which maps $x\otimes
t^s\to x\otimes(t-a)^s$ for all $x\in\lie g$, $s\ge 0$.
\end{rem}

\subl{} The following lemma is an elementary consequence of the
fact that the relations which  define the Weyl modules are graded.
\begin{lem} \hfill

\begin{enumerit}
\item The modules $W(\lambda)$ admit a $\bn$--grading induced by the
grading on $\lie g[t]$,
$$W(\lambda)=\bigoplus_{s\in\bn}W(\lambda)_s. $$
\item For all $s\ge 0$, the subspaces $W(\lambda)_s$ are
finite--dimensional $\lie g$--submodules, and we have,
$$W(\lambda)=\bigoplus_{(\mu,s)\in P\times \bn}
W(\lambda)_{\mu,s},$$  where $W(\lambda)_{\mu,s}=\{w\in
W(\lambda)_s: hw=\mu(h)w\ \forall\ h\in \lie h.\}$.
\end{enumerit}\hfill\qedsymbol
\end{lem}
 For $\mu \in P$ we set $$W(\lambda)_{\mu} =
\bigoplus_{s \in N} W(\lambda)_{\mu,s} = \{w\in W(\lambda):
hw=\mu(h)w\ \forall\ h\in \lie h.\}.$$

\begin{thm}\label{thm_uni} \cite{CPweyl} For all
$\lambda\in P^+$, the modules $W(\lambda)$ are
finite--dimensional. Moreover, any finite--dimensional $\lie
g[t]$--module $V$ generated by an element $v\in V$ satisfying the
relations
\begin{equation}\label{wquo}
\lie n^+[t]v=0,\ \ \lie ht[t]v=0,\ \ hv=\lambda(h)v,
\end{equation}
is a quotient of $W(\lambda)$.\hfill\qedsymbol
\end{thm}

\subsection{Fusion modules}

\subl{} The definition of the   fusion product of $\lie
g[t]$--modules was given in  \cite{FL} and we recall the
definition  in  the case of interest to this paper. Given
$a\in\bc$, let $ev_a:\lie g[t]\to\lie g$ be the homomorphism of
Lie algebras which  maps $x\otimes t^s$ to  $a^sx$, $x\in\lie g$,
$s\in\bn$. Let $V_a(\lambda)$ be the $\lie g[t]$--module obtained
by pulling back $V(\lambda) $ through $ev_{a}$. The following
proposition is well--known, \cite{CPnew}, \cite{FL} for instance.
\begin{prop} \label{cyclic} Let $k\in\bn_+$, $\lambda_s\in P^+$,
$a_s\in\bc$, $1\le s\le k$ and assume that $a_s\ne a_{s'}$ if
$s\ne s'$. Then, $V_{a_1}(\lambda_1)\otimes\cdots\otimes
V_{a_k}(\lambda_k)$ is an irreducible $\lie
g[t]$--module.\hfill\qedsymbol\end{prop}

\subl{} Assume  that the hypotheses of Proposition\ref{cyclic} are
satisfied and set
$$V_\boa(\blambda)=V_{a_1}(\lambda_1)\otimes\cdots\otimes
V_{a_k}(\lambda_k),\ \ \bov=v_{\lambda_1}\otimes\cdots\otimes
v_{\lambda_k}. $$ The module $V_\boa(\blambda)$ is quite clearly
not a $\bn$--graded $\lie g[t]$--module. However, the
$\bn$--grading on $\lie g [t]$ induces a $\lie g$--equivariant
filtration on $V_\boa(\blambda)$ as follows: let
$V^n_\boa(\blambda)$ be the subspace of $V_\boa(\blambda)$ spanned
by elements of the form $g\bov$, where $g\in\bu(\lie g[t])$ has
grade at most $n$.

 Set $V^{-1}_\boa(\blambda)=0$, and set
$$V_{a_1}(\lambda_1)*\cdots* V_{a_k}(\lambda_k)= \bigoplus_{n\in
\bn} V^n_\boa(\blambda),/V^{n-1}_\boa(\blambda).$$
 For $\bov' \in V_\boa(\blambda)$ let $\overline{\bov'}$ denote
its image in $V_{a_1}(\lambda_1)*\cdots* V_{a_k}(\lambda_k)$.

\begin{lem}\label{deffus} The following formula defines a
 $\bn$--graded $\lie g[t]$--module structure on
$V_{a_1}(\lambda_1)*\cdots*
 V_{a_k}(\lambda_k)$:
$$(x\otimes t^s)\overline{\bov'}=\overline{(x\otimes t^s)\bov'},$$
for all $s\in\bn$, $x\in\lie g$, $\overline{\bov'}\in
V_{a_1}(\lambda_1)*\cdots*
 V_{a_k}(\lambda_k)$.
\end{lem} \begin{pf} It suffices to prove that the action is
well--defined. But this follows,  since $(x\otimes
t^s)V^n_\boa(\blambda)\subset V^{n+s}_\boa(\blambda)$. \end{pf}
\begin{rem} The resulting $\lie g[t]$--module is called the {\em
fusion product} of the modules $V_{a_i}(\lambda_i)$, $1\le i\le
k$.\end{rem}

\subl{} The next lemma is an immediate consequence of Proposition
\ref{cyclic}.
\begin{lem} \label{fuscyclic} We have, $$V_{a_1}(\lambda_1)*\cdots*
V_{a_k}(\lambda_k)=\bu(\lie
g[t])\overline\bov.$$\hfill\qedsymbol\end{lem}

In this paper we prove a significant case of the following
conjecture.

\begin{conj}\cite{FL}\label{conj_FL} Let $k\in\bn_+$,
$\lambda_s\in P^+$, $1\le s\le k$. Assume that $a_s,b_s\in\bc$,
$1\le s\le k$ and that $a_s\ne a_{s'}$, $b_s\ne b_{s'}$ if $1\le
s\ne s'\le k$. Then, $$V_{a_1}(\lambda_1)*\cdots*
V_{a_k}(\lambda_k)\cong V_{b_1}(\lambda_1)*\cdots*
V_{b_k}(\lambda_k)$$ as $\lie g[t]$--modules.
\end{conj}\hfill\qedsymbol

  \subl{} We now  prove,
\begin{prop}\label{thm_fus}
Suppose that $k\in\bn_+$, $\lambda_i\in P^+$, $1\le i\le k$ and
that  $a_1,\cdots ,a_k\in\bc$ are distinct. The fusion product
$V_{a_1}(\lambda_1) * \dots *V_{a_k}(\lambda_k)$ is a quotient of
$W(\lambda_1 + \dots + \lambda_k)$.
\end{prop}

\begin{pf}
It suffices  in view of Lemma \ref{fuscyclic} to show that
$\overline\bov$ satisfies the relations in \eqref{wquo}. Set
$\lambda=\sum_{i=1}^k\lambda_i$. Since, $(x_{i,j}^+\otimes t^s)
\bov=0$ for all $1\le i\le j\le r$ and all $s\ge 0$, it follows
that  $(x_{i,j}^+\otimes t^s) \overline{\bov}=0$. Further,
$h\overline\bov=\lambda(h)\overline\bov$.
 Finally,
we have $$(h\otimes
t^{s+1})\bov=\left(\sum_{j=1}^k\lambda_j(h)a_j^{s+1}\right) \bov
,$$ which implies that  $(h_i\otimes t^{s+1})\overline{\bov}=0,$
for all $s\ge 0$, thus proving the proposition.

\end{pf}

\begin{cor}\label{cor_ineq}
For $\lambda=\sum_{i=1}^rm_i\omega_i$ we have $$\dim W(\lambda)\ge
\prod_{i=1}^r\binom{r+1}{i}^{m_i}.$$
\end{cor}

\begin{proof} The proof is immediate  from
Proposition~\ref{thm_fus} by taking $k=\sum_{i=1}^r m_i$,
$\lambda_s=\omega_{i_s}$, $1\le s\le k$ and we assume that each
$\omega_i$ occurs $m_i$ times.
\end{proof}
The corollary was proved originally in \cite{CPweyl} by passing to
the quantum case and the proof given there  was much more
complicated.

\subsection{Demazure modules}

 \subl{} \def \whg
{{\widehat{\lie g}}}
\def \whh {{\widehat{\lie h}}}
We now show that the Demazure modules in the highest weight
integrable representations of the affine algebra associated to
$\lie{sl}_{r+1}$ are quotients of  Weyl modules.   We begin by
recalling the definition of the affine Kac--Moody algebra. Thus,
let $\bc[t,t^{-1}]$ be the ring of Laurent polynomials in an
indeterminate $t$. Let $\whg$ be the affine Lie algebra defined by
$$\widehat{\lie g}=\lie g\otimes\bc[t,t^{-1}]\oplus\bc c\oplus \bc
d,$$ where $c$ is central and $$[x\otimes t^s, y\otimes
t^k]=[x,y]\otimes t^{s+k}+s\delta_{s+k,0}\left<x,y\right>c,\ \
[d,x\otimes t^s]=sx\otimes t^s,$$ for all $x,y\in\lie g$,
$s,k\in\bz$ and where $\left<\ ,\ \right>$ is the Killing form of
$\lie g$. Set $$\whh=\lie h\oplus\bc c\oplus\bc d, $$ and regard
$\lie h^*$ as a subspace of  $\whh^*$ by setting
$\lambda(c)=\lambda(d)=0$  for $\lambda \in \lie h^*$.

 For $0\le i\le r$,
define elements $\Lambda_i\in \whh^*$, by
$$\Lambda_i(h_j)=\delta_{i,j},\ \ \Lambda_i(d)=0,\ \
\Lambda_i(c)=1, \ \ 1\le j\le r,$$  and let $\widehat{P}^+\subset
\whh^*$ be the non--negative integer linear span of $\Lambda_i$,
$0\le i\le r$. Let $\delta\in(\widehat{\lie h})^*$ be defined by
setting $$\delta|_{\lie h} =0, \ \ \delta(d)=1,\ \ \delta(c)=0.$$
The simple roots for $\whg$ with respect to $\whh$ are $\alpha_i$,
$0\le i\le r$, where $\alpha_0=\delta-(\alpha_1+\cdots+\alpha_r)$.
Let $\widehat{Q}^+$ be the subset of $\widehat{P}^+$ spanned by
the elements $\alpha_i$, $0\le i\le r$. Let $\widehat{W}$ be the
affine Weyl group and $( , )$ be the $W$--invariant form on $\whh
^*$ obtained by requiring $(\Lambda_i,\alpha_j)=\delta_{ij} $ for
all $0\le i,j\le r$ and $(\delta,\alpha_i)=0$ for all $1\le i\le
r$. The following Lemma is well--known, see \cite{FLitt} for
instance.
\begin{lem}\label{lem_basic} Let $\mu\in P^+$. There exists $0\le
i_\mu\le r$ and $w(\mu)\in\widehat{W}$ such that
$w(\mu)\Lambda_{i_\mu}|_{\lie h}=\mu.$\hfill\qedsymbol
\end{lem}
\subl{}

Given $\Lambda=\sum_{i=0}^r m_i\Lambda_i\in\widehat{P}^+$, let
$L(\Lambda)$ be the $\whg$--module generated by an element
$v_\Lambda$ with relations: $$\lie gt[t] v_\Lambda=0,\ \quad (\lie
n^+ \otimes 1)v_\Lambda=0,\quad (h \otimes 1) v_\Lambda=\Lambda(h)
v_\Lambda, $$ $$(x_{i,i}^- \otimes 1)^{m_i+1}v_\Lambda=0,\ \
(x_{1,r}^+\otimes t^{-1})^{m_0+1}v_\Lambda=0,$$ for $1\le i\le r$.
This module is known to be irreducible and integrable (see
\cite{K}). The next proposition also can be found in \cite{K}.
\begin{prop} \label{wtL} Let $\Lambda\in \widehat{P}^+$.
\begin{enumerit}
\item  We have
$$L(\Lambda)=\bigoplus_{\Lambda'\in
\widehat{P}}L(\Lambda)_{\Lambda'},\qquad \mbox{where} \quad
L(\Lambda)_{\Lambda'}=\{v\in L(\Lambda): hv=\Lambda'(h)v,\ \
h\in\whh\}.$$ Moreover, $\dim(L(\Lambda)_{\Lambda'})<\infty$,
$\dim(L(\Lambda))_{\Lambda}=1$, and  $L(\Lambda)_{\Lambda'}=0$ if
$\Lambda-\Lambda'\notin\widehat{Q}^+$.
\item The set
\begin{equation*} \wt(L(\Lambda))=\{\Lambda'\in \widehat{P}:
L(\Lambda)_{\Lambda'}\ne 0\}\subset \Lambda-\widehat{Q}^+.
\end{equation*}
is preserved by $\widehat{W}$ and $$\dim(L(\Lambda)_{w\Lambda'} )=
\dim(L(\Lambda)_{\Lambda'} )\ \ \forall\ \ w\in\widehat{W},\ \
\Lambda'\in \widehat{P}.$$ In particular,
$\dim(L(\Lambda)_{w\Lambda} )=1$ for all $w\in\widehat{W}$.
\end{enumerit}\hfill\qedsymbol
\end{prop}

\subl{}  Given $w\in\widehat{W}$, let $v_{w\Lambda}$ be a
non--zero element in $L(\Lambda)_{w\Lambda}$ and let
$$D(w\Lambda)=\bu(\lie g[t])v_{w\Lambda}\subset L(\Lambda).$$  The
modules $D(w\Lambda)$, $w\in \widehat{W}$  are called the {\em
Demazure modules}.  Note that if  $w,w'\in \widehat{W}$ are such
that
 $w^{-1}w'\in W\subset \widehat{W}$  then $D(w\Lambda) =D(w'\Lambda).$
 It is clear from Proposition~\ref{wtL} that $D(w\Lambda)$ is
finite--dimensional for all $w\in \widehat{W}$. We can now prove,
\begin{prop}\label{demquot}
Let $\Lambda\in\widehat{P}^+$ and $w\in\widehat{W}$ be such that
$w\Lambda|_{\lie h}\in P^+$. Then, $D(w\Lambda)$ is a quotient of
$W(w\Lambda|_{\lie h})$.
\end{prop}

\begin{pf}
To prove the proposition it suffices to  show that $v_{w\Lambda}$
satisfies~\eqref{wquo}.  Suppose that $$(x_{i,j}^+\otimes
t^s)v_{w\Lambda}\ne 0,$$ for some $1\le i\le j\le r$ and $s\in
\bn$. It follows that
$w\Lambda+\alpha_{ij}+s\delta\in\wt(L(\Lambda))$. By
Proposition~\ref{wtL}~(ii), we get
$\Lambda+w^{-1}(\alpha_{i,j})+s\delta\in \wt(L(\Lambda))$ and so
$- w^{-1}\alpha_{i,j} -s\delta \in \widehat{Q}^+$. This implies in
particular, that $- w^{-1}\alpha_{i,j}\in\widehat{Q}^+$, and we
get $$(w\Lambda,\alpha_{i,j})=(\Lambda, w^{-1}\alpha_{i,j})
 \le (\Lambda, w^{-1}\alpha_{i,j}+s\delta)
\le 0$$
 which  contradicts $w\Lambda|_{\lie h}\in P^+$ unless
$(\Lambda, w^{-1}\alpha_{i,j})
 =(\Lambda, w^{-1}\alpha_{i,j}+s\delta)
= 0$. In this case, it follows  that
 the reflection $s_{w^{-1}\alpha_{i,j}+s\delta}\in \widehat W$ fixes
 $\Lambda$, and so we get
$$ \Lambda - (w^{-1}\alpha_{i,j}+s\delta) =
s_{w^{-1}\alpha_{i,j}+s\delta} ( \Lambda +
w^{-1}\alpha_{i,j}+s\delta) \in \wt(L(\Lambda)),$$ which implies
$w^{-1}\alpha_{i,j} + s\delta \in \widehat{Q}^+$,  contradicting
$- w^{-1}\alpha_{i,j} -s\delta \in \widehat{Q}^+$. Hence
$$(x_{i,j}^+\otimes t^s)v_{w\Lambda}=0,\ \ \forall\ \ 1\le i\le j\le
r,\ \ s\in \bn.$$

If $(h\otimes t^s)v_{w\Lambda}\ne 0$, then by a similar argument,
one sees that $\Lambda+s\delta\in\wt(L(\Lambda))$ which is
impossible if $s>0$  since  $\delta \in \widehat{Q}^+$. Since
$D(w\Lambda)$ is  a finite--dimensional $\lie g[t]$--module the
proposition is now immediate from \eqref{wquo}.
\end{pf}

\subl{} The following is a special case of a  result  that  can be
found in \cite{KMOTU},\cite{M},\cite{FLitt}.
\begin{thm}\label{demdim}  Assume that $\Lambda=\Lambda_{i_0}$, $0\le
i_0\le r$ and that $w\in \widehat{W}$ is such that
$w\Lambda|_{\lie h}=\sum_{i=1}^rm_i\omega_i$, where $m_i\in\bn$.
Then, we have an isomorphism of $\lie g$--modules,
$$D(w\Lambda_{i_0})\cong V(\omega_1)^{m_1}\otimes\cdots\otimes
V(\omega_k)^{m_k}.$$\hfill\qedsymbol
\end{thm}
\subsection{Main theorem}

\subl{} The next result which was conjectured in \cite{CPweyl} and
made precise in  \cite{CP2}  is the main theorem of this paper.
\begin{thm}\label{conjtrue} Let $\lambda=\sum_{i=1}^r m_i\omega_i\in
P^+$. We have, \begin{equation}\label{dimweyl} \dim (W(\lambda))=
\prod_{i=1}^r\binom{r+1}{i}^{m_i}.\end{equation}
\end{thm} We prove Theorem \ref{conjtrue} in the following
sections. We conclude this section by noting some corollaries of
the theorem.  Notice that the next corollary establishes the
Conjecture~\ref{conj_FL} for a substantial family of modules.

\begin{cor}\label{isom}
  Given $i_1,\cdots ,i_k\in\{1,\cdots ,r\}$, set
  $\lambda=\sum_{s=1}^k\omega_{i_s}$. Let $w(\lambda)\in \widehat{W}$
and
  $0\le i_\lambda\le  r$ be
  as in Lemma~\ref{lem_basic}, so that
$w(\lambda)\Lambda_{i_\lambda}|_{\lie h} =\lambda$. Let  $a_1,
\dots, a_k \in \bc$ be such that $a_i\ne a_j$,  $1\le i\ne j\le
k$. Then  we have an isomorphism of
 $\lie
g[t]$--modules $$V_{a_1}(\omega_{i_1})*\cdots
*V_{a_k}(\omega_{i_k})\cong W(\lambda)\cong
D(w(\lambda)\Lambda_{i_\lambda}).$$
\end{cor}

\begin{pf} The corollary is immediate from the Theorem~\ref{conjtrue}
together with Proposition~\ref{thm_fus}, Proposition~\ref{demquot}
and Theorem~\ref{demdim}.\end{pf}

\begin{rem} To prove the conjecture of \cite{FL} for arbitrary
irreducible representations by identifying it with a $\lie
g[t]$-module   which is obviously independent of the points, would
involve identifying a suitable quotient of the Weyl modules. A
possible approach is proposed in \cite{FKL}.
\end{rem}

\subl{} For a partition   $\xi=(\xi_1\ge\xi_2\ge\cdots \ge
\xi_{r+1})$ of non--negative integers,  set
$\lambda_\xi=\sum_{i=1}^r (\xi_i - \xi_{i+1})\omega_i$.  Given a
partition $\xi$ and regarding $\mu = \sum_{i=1}^r n_i \omega_i\in
P^+$ as the vector $(n_1, \dots, n_r)$,  let $K_{\mu, \xi}(t)$ be
the corresponding Kostka polynomial (see \cite{Mac}). Note that
$K_{\mu, \xi}(t)=0$ unless $\xi_{r+2}=0$.

 \begin{cor}\label{cor_Kostka} Retain the notation of Corollary
 \ref{isom}. Then,  $$\ch_t(\, V(\omega_{i_1})_{a_1}* \dots *
V(\omega_{i_k})_{a_k} ) =\ch_t\,  W(\lambda) = \sum_\xi
K_{\lambda, \xi^{tr}}(t) \cdot \ch_{\lie g} \, V(\lambda_{\xi}),
$$ where the sum is over all partitions $\xi$ satisfying
$\xi_{r+2}=0$,  and $\xi^{tr}$ is the transposed partition.
\end{cor}

\begin{proof} The second equality follows by   first using the
identification between the Weyl modules and Demazure modules given
in Corollary \ref{conjtrue} and  then using Theorem~5.2 in
\cite{KMOTU} which gives the formula for  $\ch_{t}(D(w\Lambda))$ .
 The first equality follows from the isomorphism
between the Weyl modules and the fusion product.
\end{proof}

\begin{rem}
 This corollary establishes a particular case of a conjecture
stated in \cite{FL}. In \cite{Ke}  another case of the conjecture
was proved, namely it was shown that the graded character of
module $$V(n_1\omega_1)*\cdots *V(n_k\omega_1), \qquad n_j\in\bn,
\quad 1\le j\le k,$$ is given by a similar formula based on Kostka
polynomials.
\end{rem}
 In the next section we shall give a {\em fermionic} character
formula for these modules.

\subl{}\label{subs_w} We conclude this section by showing that
Proposition~\ref{demquot} and Theorem \ref{conjtrue} are related
to Conjecture 1 in \cite{FLitt}. Namely, we shall prove that the
Demazure module is isomorphic as a $\lie g$--module to a cyclic
finite--dimensional module for $\whg$. For $a\in\bc^\times$,
$\lambda\in P^+$, define the module $W_a(\lambda)$ for $\whg$ as
the module with generator $w_{a,\lambda}$ and relations
\begin{equation}\label{weylgenrela} (\lie n^+\otimes
\bc[t,t^{-1}])w_{a,\lambda}=0,\ \ (h\otimes t^s)
w_{a,\lambda}=a^s\lambda(h)w_{a,\lambda},\ \ (x^-_{i,i}\otimes
1)^{m_i+1}w_{a,\lambda}=0,
\end{equation}
 for all $h\in\lie h$, $1\le i \le r$, $s\in\bz$. It was proved
in \cite{CPweyl} that as $\lie g[t]$--modules,
$$W_a(\lambda)\cong\bu(\lie g[t])w_{a,\lambda},$$ and that $$\dim
W_a(\lambda)=\dim W_b(\lambda),\ \ \forall\  a,b\in\bc^\times.$$
In fact, it was proved there, that the action of $x\otimes t^{-s}$
with $s\in\bn$ on $W_a(\lambda)$ was a linear combination of
elements from $\bu(\lie g[t])$.

Moreover if we pull back $W_a(\lambda)$ by the automorphism
$\tau_a:\lie g[t]\to \lie g[t]$ mapping  $x\otimes t^s$ to
$x\otimes (t+a)^s,$ we find that $\tau_a^*(W_a(\lambda))$ is a
quotient of $W(\lambda)$. Hence by Theorem~\ref{conjtrue} and the
lower bound for the dimension of
 $W_a(\lambda)$ obtained in \cite{CPweyl}
 we have, $$\tau_a^*(W_a(\lambda))\cong W(\lambda).$$ Let
$N(\lambda)$ be the $\lie g[t]$--submodule of $W(\lambda)$ such
that  the quotient is a Demazure module and let  $N_a(\lambda) =
\tau_{-a}^*(N(\lambda))$
 be the corresponding
submodule in $W_a(\lambda)$. It suffices to prove that
$N_a(\lambda)$ is  a $\lie g\otimes\bc[t,t^{-1}]$--submodule, but
this is now clear, since the action of $x\otimes t^{-s}$ with
$s\in\bn$ on $W_a(\lambda)$ was a linear combination of elements
from $\bu(\lie g[t])$.

\section{A basis of $W(\lambda)$}

We prove Theorem~\ref{conjtrue} by giving an explicit basis
  for $W(\lambda)$. In this section we define the proposed basis
of $W(\lambda)$ by identifying a suitable subset of the
Poincare--Birkhoff--Witt basis of $\bu(\lie n^-[t])$  and prove
that it has the desired cardinality.  We also  deduce the
fermionic version of a character  formula for $W(\lambda)$.

\subsection{The set $\cal B^r(\lambda)$}

\subl{}
 Let $\bF$ denote
 the set of all pairs
$(\ell,\bos)$ which satisfy: $$\ell\in\bn,\ \
 \bos=(\bos(1)\le \cdots \le \bos(\ell))\in\bn^\ell,$$
including the pair $(0, \emptyset)$, where $\emptyset$ is the
empty partition.
 Given $(\ell,\bos)\in\bF$ and
$1\le i\le j\le r$, let $\bx^\pm_{i,j}(\ell,\bos)$ be the element
of $\bu(\lie n^\pm[t])$ defined by,
 $$\bx^\pm_{i,j}(\ell,
\bos)=(x^\pm_{i,j}\otimes t^{\bos(1)})\cdots (x^\pm_{i,j}\otimes
t^{\bos(\ell)}),$$ if $\ell>0$ and $\bx^\pm_{i,j}(0, \emptyset)
=1$. From now on, given an element $((\ell_1,\bos_1),\cdots
,(\ell_j,\bos_j))$ of $\bF^j$, $j>0$, we let $(\ul,\us)$, be the
pair of $j$--tuples of integers,   $\ul=(\ell_1,\cdots ,\ell_j)$
and $\us=(\bos_1,\cdots ,\bos_j)$. Given $(\ul,\us)\in\bF^j$, set,
$$\bxj^\pm(\ul, \us) = \bx^\pm_{1,j} (\ell_1,
\bos_1)\cdots\bx^\pm_{j,j} (\ell_j, \bos_j). $$

\begin{defn}
 Let ${\cal B}^r$ be the  subset of $ \bu(\lie n^-[t])$
consisting of elements of the form
\begin{equation}\label{pbwform}
\underline{\bx}_1^-(\ul_1, \us_1) \dots \underline{\bx}_r^-(\ul_r,
\us_r)
\end{equation}
where $(\ul_j, \us_j) \in \bF^j$ for  $1\le j \le r$.
\hfill\qedsymbol
\end{defn}

In other words,  writing $(\ul_j, \us_j)$ as $((\ell_{1,j},
\bos_{1,j}), \dots, (\ell_{j,j}, \bos_{j,j}))$, the set ${\cal
B}^r$ consists of the elements
\begin{equation}\label{pbwform2}
\bx_{1,1}^-(\ell_{1,1},\bos_{1,1}) \bx_{1,2}^-(\ell_{1,2},
\bos_{1,2})\bx_{2,2}^-(\ell_{2,2},\bos_{2,2})
\bx_{1,3}^-(\ell_{1,3},\bos_{1,3}) \cdots
\bx_{r,r}^-(\ell_{r,r},\bos_{r,r}),
\end{equation}
in $\bu(\lie n^-[t])$ with $(\ell_{i,j},\bos_{i,j})\in \bF$, $1\le
i\le j\le r$.

\begin{prop}
The set ${\cal B}^r \subset\bu(\lie n^-[t])$ is a basis of
$\bu(\lie n^-[t])$.
\end{prop}

\begin{proof}
Fix an  ordered  basis of $\lie n^-[t]$ as follows: the elements
of the basis are  $x^-_{i,j}\otimes t^s$ for $1\le i\le j\le r$,
$s\in\bn$ and the order is given by, $x^-_{i,j}\otimes t^s$
precedes $x^-_{i',j'}\otimes t^{s'}$ when either  $j<j'$ or
$j=j'$, $i<i'$ or $j=j'$, $i=i'$, $s<s'$. The proposition follows
from the Poincare--Birkhoff--Witt theorem.
\end{proof}

\subl{}

For $m\in\bn$,  let
 \begin{equation}\label{defnbf} \bF(m)=\left\{(\ell,\bos): \ell>0,\ \
0\le \bos(i) \le m - \ell,\
 \ 1\le i\le \ell\right\}\cup\{(0,\emptyset)\}.\end{equation}
 If $m<0$ then we set $\bF(m)=\emptyset$.
For $\lambda = \sum_{i=1}^r m_i \omega_i\in P^+ $  and $1\le j\le
r$, let
\begin{equation}\label{defnbfj} \bF^j(\lambda)=\left\{(\ul,
\us)\in\bF^j: (\ell_i, \bos_i) \in \bF(m_i),\ \ 1 \le i \le
j\right\}.\end{equation}

 Given $\ul\in\bn^j$, set $\eta_j(\ul) = \sum_{i=1}^j \ell_i
\alpha_{i,j}\in Q^+$.  For any $h \in \lie h$ we have
$$[h,\underline{\bx}_j^\pm (\ul, \us)]  = \pm \eta_j(\ul)(h) \cdot
\underline{\bx}_j^\pm(\ul, \us).$$

\begin{defn}
For $\lambda\in P^+$, let ${\cal B}^r(\lambda)$ be the subset of $
{\cal B}^r$ consisting of those elements in ~\eqref{pbwform} which
satisfy  \begin{equation}\label{def1}(\ul_j, \us_j) \in
\bF^j\left(\lambda - \sum_{s=j+1}^{r} \eta_s(\ul_s)\right), \qquad
1\le j \le r.\end{equation}\hfill\qedsymbol
\end{defn}
Write $\ul_j=(\ell_{1,j},\cdots \ell_{j,j})\in\bn^j$, $1\le j\le
r$. Observe that if  $\lambda=\sum_{i=1}^rm_i\omega_i\in P^+$,
then    $$\sum_{s=j+1}^{r}
\eta_s(\ul_s)-\sum_{i=1}^j\left(\sum_{s=j+1}^r(\ell_{i,s}-\ell_{i+1,s})\right)\omega_i\in
\sum_{i=j+1}^r\bz\omega_i. $$  It is now easy to check that
\eqref{def1} is equivalent to $$(\ell_{i,j}, \bos_{i,j})
\in
\bF\left(m_i+\sum_{s=j+1}^r(\ell_{i+1,s}-\ell_{i,s} )\right),$$ or
even more explicitly to: either $\ell_{i,j}=0$, or $\ell_{i,j}>0$
and
\begin{equation}\label{adm}
\bos_{i,j}(\ell_{i,j})\le
m_i+\sum_{s=j+1}^r\ell_{i+1,s}-\sum_{s=j}^r\ell_{i,s},
\end{equation}
for all $1\le i\le j\le r$.

 Note that  for
$\lambda\in P^+$, and $\ul_r\in\bn^r$, we have
$$\cal B^{r-1}(\lambda-\eta_r(\ul))\subset\cal B^r(\lambda).$$

\begin{prop}\label{prop_count}
\begin{enumerit}
\item We have
$$\cal{B}^r(\lambda)= \bigcup_{(\ul, \us) \in \bF^r(\lambda)}
\cal{B}^{r-1}\left(\lambda-\eta_r(\ul)\right)\bxu^-(\ul, \us),$$

\item We have
$$\left|{\mathcal B}^r(\lambda)\right| =
\prod_{i=1}^r\binom{r+1}{i}^{m_i}.$$\end{enumerit}
\end{prop}

\begin{proof} Part (i) is clear.  For (ii),  recall that for
a fixed pair $\ell$, $n$ of integers, the number of elements
$\bos$ of the form, $(\bos = (\bos(1) \le \dots \le \bos(\ell) \le
n)$ is equal to
 $\binom{n+\ell}{\ell}$, and  so we have
\begin{equation}\label{eq_bin}
|\{\bos: (\ell, \bos) \in \bF(m)\}| = \binom{m}{\ell}, \qquad
|\bF(m)| = \sum_{\ell=0}^m \binom{m}{\ell}  = 2^m.
\end{equation}
We now  proceed by induction on $r$ with \eqref{eq_bin} showing
that  induction begins  at $r=1$. Using (i) and~\eqref{eq_bin},
 we get
$$|\cal{B}^r(\lambda)|=\sum_\ul |\cal B^{r-1}\left(
\lambda-\eta_r(\ul)\right)| \prod_{i=1}^r\binom{m_i}{\ell_i},$$
where the sum is over all $\underline\ell=(\ell_1,\cdots ,\ell_r)$
such that
 $0\le \ell_{i}\le m_i$ for $1\le i\le r$. This gives,
 $$|\cal{B}^r(\lambda)| =
\sum_\ul \prod_{i=1}^{r-1} \binom{r}{i}^{m_i-\ell_i+\ell_{i+1}}
\prod_{i=1}^r\binom{m_i}{\ell_i} =\prod_{ i=1}^{r}
\sum_{\ell_i=0}^{m_i} \binom{m_i}{\ell_i}
\binom{r}{i}^{m_i-\ell_i} \binom{r}{i-1}^{\ell_i} =
 \prod_{i=1}^{r} \binom{r+1}{i}^{m_i},$$
where the last equality is a consequence of the fact that
$\binom{r+1}{i} = \binom{r}{i} + \binom{r}{i-1}$ which implies
$$\binom{r+1}{i}^n = \sum_{l=0}^n \binom{n}{l} \binom{r}{i}^{n-l}
\binom {r}{i-1}^l,  $$ and hence the proposition is proved.
\end{proof}
The following corollary is immediate from Corollary
\ref{cor_ineq}.
\begin{cor}\label{ineq2}  For $\lambda\in P^+$, we have
$$\dim(W(\lambda))\ge|\cal B^r(\lambda)|.$$\end{cor}

 \subl{} The
next theorem now obviously establishes Theorem \ref{conjtrue}.

\begin{thm}\label{weylbasis} Let $\lambda=\sum_{i=1}^rm_i\omega_i\in
P^+$. The set $\{bw_\lambda:b\in \cal B^r(\lambda)\}$ is a basis
of $W(\lambda)$.
\end{thm}
 \begin{rem} The theorem was proved in \cite{CPweyl} when $r=1$.
In the rest of the paper we shall use that case to prove Theorem
\ref{weylbasis} when $r\ge 2$.\end{rem}

Again, we postpone the proof of the theorem and deduce some
corollaries.
\begin{cor}\label{cor_fd} The elements
$(x_{1,1}^-)^{\ell_{11}}(x_{1,2}^-)^{\ell_{1,2}}
(x_{2,2}^-)^{\ell_{2,2}} (x_{1,3}^-)^{\ell_{1,3}}
(x_{2,3}^-)^{\ell_{2,3}} \cdots
(x_{r,r}^-)^{\ell_{r,r}}v_\lambda$ where $\ell_{i,j}\in\bn ^+$
satisfy
\begin{equation}\label{admf}
m_i+\sum_{s=j+1}^r\ell_{i+1,s}-\sum_{s=j}^r\ell_{i,s} \ge 0
\end{equation}
are a basis $V(\lambda)$.\end{cor}
\begin{pf}
It suffices to observe that  the subspace of $W(\lambda)$ with
$\bn$--grading zero is isomorphic to $V(\lambda)$.
\end{pf}

\subl{} Using this  basis we get a fermionic formula for the
character of $W(\lambda)$.

  Given an integer $n\in\bn $, set
$$[n]_t=\frac {1-t^n}{1-t}\in\bn[t], \qquad [n]_t! = [1]_t \dots
[n]_t, \qquad \bin{n}{s}_t = \frac{[n]_t!}{[s]_t! [n-s]_t!}.$$
Clearly all these elements are in $\bn[t]$.

\begin{prop}  We have,
\begin{equation}\label{eq_ferm}
\ch_t\,  W(\lambda) =
\sum_{{(\ell_{i,j})\in\bn^{r(r+1)/2}}\atop{1\le i\le j\le r}}
\left(\prod_{1\le i \le j \le r} \left[{ m_i +
\sum\limits_{s=j+1}^r \ell_{i+1,s} - \sum\limits_{s=j+1}^r
\ell_{i,s}}\atop{\ell_{i,j}} \right]_t\right)
e\left(\lambda-\sum_{1\le i\le j\le
r}\ell_{i,j}\alpha_{i,j}\right).
\end{equation}
\end{prop}

\begin{proof} By Theorem \ref{weylbasis}, it follows that
$$\ch_t(W(\lambda))=\sum_{b\in\cal
B^r(\lambda)}e(\mu(b))t^{s(b)},$$ where $\mu(b)\in P$,
$s(b)\in\bn$ are defined as follows. If
\begin{equation*} b=\bx_{1,1}^-(\ell_{1,1},\bos_{1,1})
\bx_{1,2}^-(\ell_{1,2},
\bos_{1,2})\bx_{2,2}^-(\ell_{2,2},\bos_{2,2})
\bx_{1,3}^-(\ell_{1,3},\bos_{1,3}) \cdots
\bx_{r,r}^-(\ell_{r,r},\bos_{r,r}),
\end{equation*} for some  $(\ell_{i.j},\bos_{i,j})\in\bF$, $1\le i \le
j \le r$,  then $$s(b)= \sum_{1\le i\le j\le r
}\sum_{m=1}^{\ell_{i,j}}\bos_{i,j}(m),\qquad
\mu(b)=\lambda-\sum_{1\le i\le j\le r}\ell_{i,j}\alpha_{i,j}.$$
Hence we get $$\ch_t\,  W(\lambda) =
\sum_{{(\ell_{i,j})\in\bn^{r(r+1)/2}}\atop{1\le i\le j\le r}}
e\left(\lambda-\sum_{1\le i\le j\le
r}\ell_{i,j}\alpha_{i,j}\right) \prod_{1\le i\le j\le r\ } \sum_{\
\bos : \, (\ell_{i,j},\bos)\in\bF(m_{i,j})}
 t^{\bos(1) + \dots + \bos(\ell)},$$
where $m_{i,j} = m_i + \sum\limits_{s=j+1}^r \ell_{i+1,s} -
\sum\limits_{s=j+1}^r \ell_{i,s}$. But now, the result follows by
observing that for a fixed pair of integers  $\ell, m \in \bn$  we
have $$\sum_{\bos : \, (\ell,\bos)\in\bF(m)} t^{\bos(1) + \dots +
\bos(\ell)} = \sum_{0 \le \bos(1) \le  \dots \le \bos(\ell) \le
m-\ell} t^{\bos(1) + \dots + \bos(\ell)} = \bin{m}{\ell}_t.$$

\end{proof}

\begin{rem} It  follows from  \cite[Proposition~5.3 ]{HKKOTY} that
the right hand side of~\eqref{eq_ferm} is equal to $\sum
K_{\lambda, \xi^{tr}}(t) \cdot \ch_{\lie g} \, V(\lambda_{\xi})$.
Thus, the preceding proposition gives an alternate proof of
\cite[Theorem 5.2]{KMOTU}.\end{rem}

\subsection{A Gelfand--Tsetlin type  filtration for $W(\lambda)$}
 In this   section we construct a filtration indexed by $\bn^{2r}$
of $W(\lambda)$ and show that Theorem~\ref{weylbasis} follows from
Proposition~\ref{assgraded} which studies  the associated graded
space of the filtration. As a consequence, we  see that  the
associated graded module is isomorphic to a direct sum of Weyl
modules for $\lie \lie{sl}_r[t] \subset \lie g[t]$.

\subl{} Let $\lie g_{r-1}$ be the subalgebra of $\lie g$
isomorphic to $\lie{sl}_r$ which is  generated by the elements
$x_{i,j}^\pm$ $1\le i\le j\le r-1$.
  Set $$\lie n^\pm_{r-1}=\bigoplus_{1\le i\le j\le r-1}\bc
x_{i,j}^\pm,\qquad  \lie u_r^\pm=\bigoplus_{i=1}^r\bc
x_{i,r}^\pm.$$
 Clearly the
elements $\{\bxu^\pm(\ul,\us):(\ul,\us)\in\bF^r\} $ are a basis
for $\bu(\lie u_r^\pm[t])$, and we have by using the PBW--theorem
that $$W(\lambda)=\sum_{(\ul,\us)\in \bF^r}\bu(\lie
g_{r-1}[t])\bxu^-(\ul,\us).$$

\subl{}
 Let $\bi^r=\bn^r\times \bn^r$. Introduce an order on
$\bi^r = \{(\ul, \ud)\}$ in the following way. We say that
\begin{eqnarray*}
\ul \eql \ul' & \mbox{if} & \ell_1 = \ell_1',\dots, \ell_{s-1} =
\ell_{s-1}',\ \ell_s < \ell_s',\\ \ud \eqd \ud' & \mbox{if} & d_r=
d_r',\dots, d_{s+1} = d_{s+1}',\ d_s
> d_s'\end{eqnarray*} for some $1\le s\le r$. Finally we say that $$
(\ul,\ud)> (\ul', \ud') \ \ \mbox{if} \ \  \ul \eql \ul' \
\mbox{or} \  \ \  \ul = \ul', \ \  \ud \eqd \ud'.$$

For, $(\ell, \bos) \in \bF$ set
$|\bos|=\sum_{p=1}^{\ell}\bos(p)$, if $\ell>0$ and
$|\emptyset|=0$.
Given $(\ul, \us) \in \bF^r$ let $|\us| = (|\bos_1|, \dots,
|\bos_r|)$. Note that $(\ul,|\us|)\in\bi^r$.
 Given $\boi\in\bi^r$, define the $\lie g_{r-1}[t]$--modules,
\begin{eqnarray*}&W(\lambda)^{\ge \boi}
&=\sum_{\{(\ul,\us)\in\bF^r: (\ul,|\us|)\ge \boi\}}\bu(\lie
g_{r-1}[t])\bxu^-(\ul,\us)w_\lambda,\\ & W(\lambda)^{> \boi}
&=\sum_{\{(\ul,\us)\in\bF^r: (\ul,|\us|)> \boi\}}\bu(\lie
g_{r-1}[t])\bxu^-(\ul,\us)w_\lambda ,\\ &
\gr(W(\lambda))&=\bigoplus_{\boi\in\bi^r}W(\lambda)^{\ge\boi}/W(\lambda)^{>\boi}.\end{eqnarray*}

Let $$\gr^\boi:
W(\lambda)^{\ge\boi}\to W(\lambda)^{\ge\boi}/W(\lambda)^{>\boi}$$
the canonical projection which is clearly a map of $\lie
g_{r-1}[t]$--modules. The next result indicates the usefulness of
this construction.

\begin{prop}\label{prop_grsp} For all $\lambda\in P^+$, we have an
isomorphism of vector spaces, $$W(\lambda)\cong\gr(W(\lambda)).$$
More precisely,
 let $\{v_1,\cdots ,v_N\}\in W(\lambda)$ be a set of vectors such
that $v_j\in W(\lambda)^{\boi_j}$ for some $\boi_j\in\bi^r$, $1\le
j\le N$. Suppose that $\{\gr^{\boi_1}(v_1),\cdots
,\gr^{\boi_N}(v_N)\}$ span
 $\gr(W(\lambda))$. Then $\{v_1,\cdots ,v_N\}$  span
 $W(\lambda)$.\end{prop}

\begin{pf}
Let $\bi^r(\lambda)$ be the subset of $\bi^r$ defined by,
$$\bi^r(\lambda) = \left\{\boi \in \bi^r:
\sum_{\{(\ul,\us)\in\bF^r: (\ul,|\us|)= \boi\}} \bu(\lie
g_{r-1}[t])\bxu^-(\ul,\us)w_\lambda \ne 0\right\}.$$ We first show
that $\bi^r(\lambda)$ is finite. For this, it suffices to prove
that $\bx_r^-(\ul,\us)w_\lambda=0$ for all but finitely many pairs
$(\ul,\us)\in\bF^r$. Note that
\begin{equation*}
\bx_r^-(\ul,\us)w_\lambda\in W(\lambda)_{\mu_\ul,k_\us},\qquad
\mu_\ul=\lambda-\eta_r(\ul), \quad k_\us
=\sum_{i=1}^r|\bos_i|.
\end{equation*}
Since $W(\lambda)$ is a finite--dimensional $\lie g[t]$--module,
we can choose a finite set of pairs $(\mu_1, k_1), \dots, (\mu_M,
k_M)$ such that $W(\lambda)_{\mu,k}=0$ unless $\mu = \mu_j$ and
$k= k_j$ for some $1\le j \le M$. Then it is enough  to show that
the sets $$\bs_j=\{(\ul,\us)\in\bF^r:\mu_\ul=\mu_j,\ \
k_\us=k_j\},\ \ 1\le j\le M,$$ are finite. Since the elements
$\alpha_{i,r}\in Q^+$, $1\le i\le r$ are linearly  independent, it
follows immediately that for $1\le j\le M$, either $\bs_j$ is
empty (in which case there is nothing to prove), or there exists
unique $\ul(j)\in\bn^r$, such that $\mu_{\ul(j)} = \mu_j$, so
$$\bs_j=\{(\ul(j),\us)\in\bF^r: |\bos_1|+ \dots + |\bos_r|
=k_j\}.$$ Since  each  $\bos_i$ is a tuple of non-negative
integers, it follows now that $\bs_j$ is finite for all $1\le j\le
M$.

Let  $\boi_1 <\boi_2,< \cdots <\boi_R$ be an ordered enumeration
of the elements of $\bi^r(\lambda)$. Let $\boi_{R+1}$ be any
element of $\bi^r$ greater than $\boi_R$. Then we have
$$W(\lambda)^{\ge \boi_1} = W(\lambda), \qquad W(\lambda)^{\ge
\boi_{R+1}} = 0, \qquad W(\lambda)^{>\boi_k} = W(\lambda)^{\ge
\boi_{k+1}}, \quad 1\le k \le R.$$

Thus, we get $$\gr (W(\lambda)) = \bigoplus_{\boi \in
\bi^r(\lambda)} W(\lambda)^{\ge \boi}/W(\lambda)^{> \boi} =
\bigoplus_{k=1}^r W(\lambda)^{\ge \boi_{k}}/W(\lambda)^{\ge
\boi_{k+1}},$$ which is the associated graded space of a usual
increasing filtration $$0\subset W(\lambda)^{\ge
\boi_R}\subset\cdots\subset W(\lambda)^{\ge \boi_1}=W(\lambda),$$
and for such a filtration the statement of the proposition is
standard.

\end{pf}

\subl{} We shall deduce Theorem \ref{weylbasis} from the next
proposition.

\begin{prop}{\label{assgraded}}  Let $r\ge
2$.
\begin{enumerit}
\item  For  $(\ul,\us)\in\bF^r$,
 there exists a map of $\lie g_{r-1}[t]$--modules
$\psi^{\ul,\us}: W(\lambda-\eta_r(\ul))\to\gr(W(\lambda))$ given
by extending  $$\psi^{\ul,\us}(w_{\lambda-\eta_r(\ul)})
=\gr^{(\ul, |\us|)}(\bxu^-(\ul,\us)w_\lambda).$$ Moreover,
$$\gr(W(\lambda))=\bigoplus_{(\ul,\us)\in\bF^r}{\text{Im}}(\psi^{\ul,\us}).$$

\item  The images of $\psi^{\ul,\us}$  for
 $(\ul,\us)\in\bF^r(\lambda)$
span $\gr(W(\lambda))$.

\end{enumerit}
\end{prop}

\subl{}\label{sec_mtproof}{\em Proof of Theorem \ref{weylbasis}}.
The theorem is proved by induction on $r$. Note that induction
starts at $r=1$ by \cite[Theorem , Section 6]{CPweyl}.  Assume the
result for $r-1$. Using Proposition~\ref{assgraded} and the
induction hypothesis we see that $\gr(W(\lambda))$ is spanned by
the  sets
$$\gr^{(\ul,|\us|)}(\cal{B}^{r-1}(\lambda-\eta_r(\ul))\bxu^-(\ul,\us)
w_\lambda), \qquad \mbox{where} \quad (\ul, \us)\in\bF(\lambda).$$
Proposition~\ref{prop_grsp} together  with
Proposition~\ref{prop_count}~(i) now shows that
$\cal{B}^r(\lambda)$ spans $W(\lambda)$ and hence $\dim
W(\lambda)\le\ |\cal{B}^r(\lambda)|$. Since Corollary \ref{ineq2}
established the reverse inequality, it follows that $\dim
W(\lambda)= |\cal{B}^r(\lambda)|$. Hence $B^r(\lambda)$ is a basis
of $W(\lambda)$ and the theorem is proved. \hfill\qedsymbol

\subl{} We conclude this section with a final corollary of
Theorem~\ref{weylbasis} and Proposition~\ref{assgraded}.

\begin{cor} For $r\ge 2$ and  $\lambda = \sum
m_i \omega_i$, we have an isomorphism of $\lie
g_{r-1}[t]$--modules
 $$\gr(W(\lambda)) \cong\bigoplus_{{\ul\in\bn^r}\atop{\ell_i \le m_i}}
m_{\ul,\lambda}W(\lambda-\eta_r(\ul)),\qquad \mbox{where} \quad
 m_{\ul, \lambda}=\binom{m_1}{\ell_1} \cdots
\binom{m_r}{\ell_r}.  $$  In particular, the module $W(\lambda)$
admits a filtration of $\lie g_{r-1}[t]$--modules such that the
composition factors are Weyl modules for $\lie g_{r-1}[t]$.
\hfill\qedsymbol
\end{cor}

\section{Proof of Proposition~\ref{assgraded}}
The proof of the proposition is quite complicated and requires a
fairly detailed analysis of the structure of $\bu(\lie g[t])$ and
of the canonical projection $\bpr: \bu(\lie g[t])\to\bu(\lie
n^-[t])$ defined below.  One has to study in some detail the
commutation relations in the algebras and their behavior with
respect to the ordered set $\bi^r$. Section~3.1 is a  collection
of  elementary properties of $\bpr$, Section~3.2 analyses the
$\lie g_{r-1}[t]$--module structure of $W(\lambda)$ and concludes
by proving  Proposition~\ref{assgraded}~(i).  Section~3.3 studies
the finite--dimensional irreducible $\lie g$--module $V(\lambda)$.
The point of this section is to obtain a suitable spanning set for
this module and to use it later in the study of $W(\lambda)$ which
we recall can be written as $W(\lambda)=\bu(\lie
gt[t])V(\lambda)$. Section~3.4 is devoted to the study of
$\bu(\lie \lie{sl}_2[t])$, more precisely  it describes the subspace of
$\bu(\lie n^-[t])$  annihilating $w_\lambda$. This is then used in
the  general case, to divide $\bu(\lie n^-[t])$ into two different
subspaces whose elements have different behavior when applied to
$w_\lambda$.
 After some further technical
results,  including the crucial lemma in Section~3.5, which
calculates $\bpr$ modulo higher terms in the Gelfand-Tsetlin
filtration, the proof of  Proposition~\ref{assgraded}~(ii) is
completed in Section~3.6.

\subsection{ The projection $\bpr$.}

\subl{}  Let $\lie b^+$ be the subalgebra $\lie h\oplus\lie n^+$
of $\lie g$ and note that $\lie g = \lie b^+ \oplus \lie n^-.$ Let
$\bpr: \bu(\lie g[t])\to\bu(\lie n^-[t])$ be the projection
corresponding to the vector space decomposition $$\bu(\lie
g[t])=\bu(\lie n^-[t])\oplus \bu(\lie g[t])(\lie b^+[t]),$$ given
by the Poincare--Birkhoff--Witt theorem. Clearly $\bpr$ is a
$\bn$--graded linear map. The next result collects some properties
of $\bpr$ which are immediate from the definition.
\begin{prop}\label{proppr}
\begin{enumerit}\item  For all $g^-\in\bu(\lie
n^-[t])$ and $x\in\bu(\lie g[t])$,  we have  $
\bpr(g^-x)=g^-\bpr(x)$.
\item  For all
$g^+\in\bu(\lie b^+[t])(\lie b^+[t])$, $x\in\bu(\lie g[t])$, we
have $\bpr(xg^+)=0$ and hence $\bpr(g^+x)=\bpr([g^+,x])$.
\item  For all $g_1, \, g_2 \in  \bu(\lie g[t])$ we have $\bpr(g_1 g_2) = \bpr(g_1
\bpr(g_2))$.
\end{enumerit}
\end{prop}\hfill\qedsymbol

Note that $\lie n^-[t] \oplus \lie b^+ t[t]$ is a subalgebra of
$\lie g[t]$.
\begin{lem} \label{prl} Let $\lambda\in P^+$. Then for all
$x \in \bu(\lie n^-[t] \oplus \lie b^+ t[t])$ we have $xw_\lambda
= \bpr(x) w_\lambda$ in $W(\lambda)$.
\end{lem}

\begin{proof} The proof is immediate from the observation that
$(\lie n^+ [t]\oplus\lie ht[t]) w_\lambda =0$.
\end{proof}

\subl{}

Let
\begin{eqnarray*}
&\bu(\lie n^-[t])^{\boi} & = \sum_{\{(\ul,\us)\in\bF^r:
(\ul,|\us|) = \boi\}}\bu(\lie n^-_{r-1}[t])\bxu^-(\ul,\us)\\
&\bu(\lie n^-[t])^{\ge \boi} &=\sum_{\{(\ul,\us)\in\bF^r:
(\ul,|\us|)\ge \boi\}}\bu(\lie n^-_{r-1}[t])\bxu^-(\ul,\us),\\ &
\bu(\lie n^-[t])^{> \boi} &=\sum_{\{(\ul,\us)\in\bF^r:
(\ul,|\us|)> \boi\}}\bu(\lie n^-_{r-1}[t])\bxu^-(\ul,\us).
\end{eqnarray*}

Clearly we have $$W(\lambda)^{\ge \boi} \supset \bu(\lie
n^-[t])^{\ge \boi} w_\lambda,\ \ \ W(\lambda)^{> \boi}\supset
\bu(\lie n^-[t])^{> \boi} w_\lambda.$$  We will see in
Proposition~\ref{act2} that  the reverse inclusions hold as well.

\subl{}

The algebra $ \bu(\lie g[t])$ is  a sum of eigenspaces for the
adjoint action of $\lie h$ on $\lie g[t]$ and the  decomposition
is compatible with the $\bn$--grading on $\bu(\lie g[t])$. More
precisely, we can write  $$\bu(\lie g[t])=\bigoplus_{\mu\in\lie
h^*} \bu(\lie g[t])_\mu =\bigoplus_{\mu\in\lie
h^*,s\in\bn}\bu(\lie g[t])_{\mu,s}.$$ The subspaces $\bu(\lie
n^\pm[t])_{\mu,s}$ are defined analogously. Note that $\bu(\lie
g[t])_{\mu,s}=0$ if $\mu\notin Q$ and that $\bu(\lie
n^\pm[t])_{\pm\mu,s}=0$ if $ \mu\notin Q^+$. Again, the next lemma
collects together some elementary properties of these
decompositions.

\begin{lem}\label{lem_wt}
\begin{enumerit}
\item For $\lambda\in P^+$, $\mu\in P$ and $s\in\bn$, we have
$W(\lambda)_{\mu,s} = \bu(\lie n^-[t])_{\mu - \lambda,s}
w_\lambda$.
\item Let  $\eta,\eta' \in Q$, $s,s' \in \bn$. We have
$\bu(\lie g[t])_{\eta,s} \bu(\lie g[t])_{\eta',s'} \subset\bu(\lie
g[t])_{\eta+\eta',s+s'}$. Moreover, if $g_1, g_2 \in \bu(\lie
n^-[t])$, then
 $$g_1\in\bu(\lie n^-[t])_{\eta,s},\ \   g_1g_2 \in
\bu(\lie n^-[t])_{\eta',s'}\quad  \mbox{implies} \quad  g_2 \in
\bu(\lie n^-[t])_{\eta'-\eta,s'-s}.$$
\item For all $\eta\in Q$, $s \in \bn$ we have
$\bu(\lie n^-[t])_{\eta,s} \subset \bu(\lie g[t])_{\eta,s}$ and
$\bpr(\bu(\lie g[t])_{\eta,s})\subset \bu(\lie n^-[t])_{\eta,s}$.
\item For all $1\le i\le j\le r$, and $(\ell,\bos)\in\bF$, we have,
$\bx_{i,j}^\pm(\ell, \bos) \in \bu(\lie g^-[t])_{\pm \ell
\alpha_{i,j}, |\bos|}$.
\end{enumerit}
\end{lem}\hfill\qedsymbol

\subl{} \ Let $\bF_+$ denotes the subset of $\bF$ consisting of
pairs $(\ell,\bos)\in\bF$ such that  $\bos(i)>0$ for $1 \le i \le
\ell$, together with   the pair $(0,\emptyset)$. The next
proposition is used repeatedly in this section.

\begin{prop}\label{prop_prg} Let $1\le i\le r$,
 $(\ell^+, \bos^+) \in \bF_+$, $(\ell, \bos) \in \bF$. The
element
$\bpr\left(\bx^+_{i,r}(\ell^+,\bos^+)\bx^-_{i,r}(\ell,\bos)\right)$
is a linear combination of elements of the form
$\bx^-_{i,r}(\ell-\ell^+,\bos')$, where $\bos'$ satisfies the
condition  $|\bos'| = |\bos|+ |\bos^+|$.
\end{prop}
\begin{proof}   Note that the
subspace $\lie s_i[t]$  of $\lie g[t]$ spanned by the elements
$x_{i,r}^\pm\otimes t^k $, $h_i\otimes t^k$, $k\in\bn$ is a
subalgebra of $\lie g[t]$ which is  isomorphic to $\lie{sl}_2[t]$.
Further,  $\bpr$ maps $\bu(\lie s_i[t])$ onto $\bu(\lie
n^-_{i,r}[t])$, where $\lie n_{i,r}^-[t]$ is the subalgebra of
$\lie n^-[t]$ spanned by elements of the form $x^-_{i,r}\otimes
t^k$, $k\in\bn$. This proves that
$\bpr\left(\bx^+_{i,r}(\ell^+,\bos^+)\bx^-_{i,r}(\ell,\bos)\right)$
is a linear combination of elements of the form
$\bx^-_{i,r}(\ell',\bos')$ for some $(\ell',\bos')\in\bF$. But
now, we see that since
$\bx^+_{i,r}(\ell^+,\bos^+)\bx^-_{i,r}(\ell,\bos) \in \bu(\lie
g[t])_{(\ell^+-\ell^-)\alpha_{i,r},|\bos|+ |\bos^+|}$ the same
holds for
$\bpr\left(\bx^+_{i,r}(\ell^+,\bos^+)\bx^-_{i,r}(\ell,\bos)\right)$.
 This means that we can assume that  $\ell' =\ell^+-\ell$ and $|\bos'|=|\bos|+
|\bos^+|$.

\end{proof}

\subsection{Action of $\lie g_{r-1}[t]$ on $W(\lambda)$}
\subl{}\
 Let $\boe_i \in
\bn^r$, $1\le i \le r$, denote the standard basis vectors.
 \begin{lem}\label{act1} Let $(\ul,\us)\in \bF^r$, $g \in \lie b^+_{r-1}[t]$.
\begin{enumerit}
\item
The   commutator $[g, \bxu^-(\ul,\us)]$ is a linear combination of
terms of the form $\bxu^-(\ul',\us')$ with $(\ul',\us') \ge
(\ul,\us)$. Moreover, if  $g \in \lie n ^+_{r-1}[t]\oplus\lie
h_{r-1}t[t]$ then we have $(\ul',\us') > (\ul,\us)$.
\item Let $\lambda\in P^+$. The element $g\bxu^-(\ul,\us) w_\lambda\in W(\lambda)$
 is a linear combination of terms of the form
$\bxu^-(\ul',\us')w_\lambda$ with $(\ul',\us') \ge (\ul,\us)$.
Moreover, if  $g \in \lie n ^+_{r-1}[t]\oplus\lie h_{r-1}t[t]$
then we have $(\ul',\us') > (\ul,\us)$.
\end{enumerit}
\end{lem}

\begin{proof}
For (i) it is enough to consider the case $g = h \otimes t^s$,
where $h \in \lie h_{r-1}$ and $s \in \bn$, and the case $g =
x^+_{i,j} \otimes t^s$, where $1\le i\le j<r$ and $s \in \bn$.
Recall the adjoint action of $\lie g[t]$ on $\bu(\lie g[t])$ is a
derivation, i.e
\begin{equation}\label{leibn}
[y,x_1\cdots x_s]=\sum_{j=1}^s x_1\cdots
x_{j-1}[y,x_j]x_{j+1}\cdots x_s,
\end{equation}
for all $y,x_1,\cdots ,x_s\in\lie g$. Now, since $$[h\otimes
t^s,x^-_{k,r}\otimes t^m]=- \alpha_{k,r}(h)x^-_{k,r}\otimes
t^{m+s}$$ for all $h\in\lie h$, $m,s\in\bn$, $1\le k\le r$, we see
that  $\left[h \otimes t^s, \bxu^-(\ul,\us) \right]$ is a linear
combination of terms $\bxu^-(\ul',\us')$ with $\ul' = \ul$,
$|\us'| = |\us|+s\boe_j$, $1\le j \le r$  and hence
$(\ul',|\us'|)\ge (\ul,|\us|)$ with the inequality being strict if
$s>0$.

Next for $1\le i\le j<r$ we have $[x_{i,j}^+\otimes t^s,
x^-_{k,r}\otimes t^m]= - \delta_{k,i}x^-_{j+1,r}\otimes t^{m+s}$,
and it follows again that $\left[x^+_{i,j} \otimes t^s,
\bxu^-(\ul,\us) \right]$ is a linear combination of terms
$\bxu^-(\ul',\us')$ with $$\ul' = \ul+\boe_{j+1} - \boe_i, \qquad
|\us'| =  |\us| +(m+s)\boe_{j+1} - m\boe_i.$$ Further,
$(\ul',\us')
> (\ul,\us)$.

For  (ii) note that $$g\bxu^-(\ul,\us) w_\lambda = \bxu^-(\ul,\us)
g  w_\lambda + [g,\bxu^-(\ul,\us)] w_\lambda.$$ The first term
 is proportional to $w_\lambda$ and is zero if  $g \in
\lie n ^+_{r-1}[t]\oplus\lie h_{r-1}t[t]$ and hence is of the
desired form. Using part (i), we see that  the second term is also
in the correct form and the lemma is proved.
\end{proof}

\subl{} \begin{prop}\label{act2} We have $$W(\lambda)^{\ge \boi} =
\bu(\lie n^-[t])^{\ge \boi} w_\lambda,\ \ \ W(\lambda)^{> \boi} =
\bu(\lie n^-[t])^{> \boi} w_\lambda.$$\end{prop}

\begin{proof}
It is enough to show that $\bu (\lie g_{r-1} [t])  \bxu^-(\ul,\us)
w_\lambda$ is contained in $\bu(\lie n^-[t])^{\ge (\ul,\us)}
w_\lambda$.  Since, $$\bu (\lie g_{r-1} [t]) =\bu (\lie n^-_{r-1}
[t])\bu (\lie b_{r-1}^+ [t]),$$ it suffices to prove that
 $$\bu (\lie b_{r-1}^+ [t])
\bxu^-(\ul,\us) w_\lambda\subset\bu(\lie n^-[t])^{\ge (\ul,\us)}
w_\lambda.$$
 But this follows from Lemma \ref{act1}~(ii).
\end{proof}

\subl{} {\em Proof of Proposition ~\ref{assgraded}~(i)}.  Let
$(\ul,\us)\in\bF^r$.  Using Lemma~\ref{act1}~(ii), we see that for
all $h\in\lie h$, $k\in\bn$, $s\in\bn_+$ and $1\le i\le j\le r$,
we have $$ (x_{i,j}^+\otimes t^k)\cdot \gr^{(\ul,
|\us|)}(\bxu^-(\ul,\us)w_\lambda)=0,\qquad (h\otimes t^s)\cdot
\gr^{(\ul, |\us|)}(\bxu^-(\ul,\us)w_\lambda)=0,$$ $$h\cdot
\gr^{(\ul, |\us|)}(\bxu^-(\ul,\us)w_\lambda)=
\left(\lambda(h)-\sum_{j=1}^r\ell_j
\alpha_{j,r}(h)\right)\gr^{(\ul, |\us|)}(\bxu^-(\ul,\us)w_\lambda)
=(\lambda-\eta_r(\ul))(h)\gr^{(\ul,
|\us|)}(\bxu^-(\ul,\us)w_\lambda).$$
This means that the
element $\gr^{(\ul,
|\us|)}(\bxu^-(\ul,\us)w_\lambda)\in\gr(W(\lambda))$ satisfies the
relations in~\eqref{wquo} with the weight $ \lambda-\eta_r(\ul)$.
Theorem~\ref{thm_uni} now implies the existence of the map
$\psi^{\ul,\us}$. Since,
 $${\rm Im} \left(\psi^{\ul,\us}\right) = \bu(\lie
g_{r-1}[t])\cdot\gr^{(\ul, |\us|)}(\bxu^-(\ul,\us)w_\lambda) =
\gr^{(\ul, |\us|)}(\bu(\lie
g_{r-1}[t])\bxu^-(\ul,\us)w_\lambda),$$ the second statement
follows.

\subsection{The case of  $V(\lambda)$}
It is necessary in this section alone, to work with the Lie
algebra $\lie{gl}_{r+1}=\lie g\oplus\bc c$.

\subl{} It is well--known that the finite--dimensional irreducible
representations of $\lie{gl}_{r+1}$ are parameterized by partitions $\xi
= (\xi_1 \ge \dots \ge \xi_{r} \ge \xi_{r+1} \ge 0)$ of
non--negative integers. Let $V(\xi)$ denote the corresponding
representation, which is generated as a $\lie{gl}_{r+1}$--module by an
element $v_\xi$ with the relations:
\begin{equation}\label{defrel2}
x_{i,j}^+v_\xi=0,\ \ H_i v_\lambda= \xi_i v_\xi,\ \
(x_{i,i}^-)^{\xi_i-\xi_{i+1}+1} v_\xi = 0,
\end{equation}
for all $1\le i\le j\le r$. Note that $V(\xi)=\bu(\lie n^-)v_\xi$.
Moreover, we have  $V(\xi)\cong V(\lambda_\xi)$ as $\lie
g$--modules, where $$\lambda_\xi = \sum_{i=1}^{r} (\xi_i -
\xi_{i+1}) \omega_i\in P^+.$$
 Conversely, given $\lambda=\sum_im_i\omega_i\in P^+$, we have
$V(\lambda)\cong V(\xi^\lambda)$ where
$\xi^\lambda=((\xi_1^\lambda \ge \dots \ge \xi_{r}^\lambda \ge
\xi_{r+1}^\lambda \ge 0)$ is defined by
$\xi_i^\lambda=\sum_{j=i}^r m_j$ and $\xi_{r+1}^\lambda=0$.

\begin{thm}\cite{GT}\label{thm_GT}
 As a module for $\lie{gl}_{r} \subset \lie{gl}_{r+1}$ we have $$
V(\xi)=\bigoplus_\eta  V(\eta),$$ where the sum is over all
partitions $\eta= (\eta_1\ge \dots \ge \eta_{r})$ of non--negative
integers satisfying,  ${\xi_i\ge \eta_i \ge \xi_{i+1}}$ for all
$1\le i\le r$.\hfill\qedsymbol
\end{thm}
The inductive construction of the Gelfand--Tsetlin basis for the
$\lie{gl}_{r+1}$--module $V(\xi)$  (see \cite{GT}) is based on this
theorem and it motivated our definition of the Gelfand--Tsetlin
filtration of the $\lie g[t]$--modules $W(\lambda)$.

 \subl{} The next two results
are probably well--known, but we include a proof here, since it is
basic for the results of this paper.

\begin{lem}\label{lem_fd} Let $\lambda = \sum m_i \omega_i$. Then we have,
\begin{equation} V(\lambda)=\sum_{{k_1, \dots, k_r \in \bn}\atop{k_1 \le m_1}}
\bu(\lie n^-_{r-1})(x^-_{1,r})^{k_1}\cdots
(x^-_{r,r})^{k_r}v_\lambda.\end{equation}
\end{lem}
\begin{pf}
 To prove the lemma, it suffices to establish it for the $\lie{gl}_{r+1}$--module
 $V(\xi^\lambda)$,
i.e to prove that that
  $$V(\xi^\lambda)=\sum_{{k_1, \dots, k_r \in
\bn}\atop{k_1 \le m_1}} \bu(\lie n^-_{r-1})(x^-_{1,r})^{k_1}\cdots
(x^-_{r,r})^{k_r}v_{\xi^\lambda}.$$

 Using Theorem~\ref{thm_GT} we see that  $V(\xi^\lambda)$ is
spanned by  the sets $\bu(\lie n^-_{r-1}) v_{\eta}$, where
$\eta=(\eta_1\ge\cdots \ge \eta_r)$ are non--negative integers
$\xi_i^\lambda \ge \eta_i \ge \xi_{i+1}^\lambda$, $1\le i \le r$,
and  $v_{\eta}\in V(\xi^\lambda)$ satisfies
\begin{equation} \label{xir} x_{i,j}^+ v_{\eta}=0, \qquad
H_m v_{\eta}=\eta_m v_{\eta}.
\end{equation}

 Writing $v_{\eta}$ as a $\bu(\lie n^-_{r-1})$--linear combination of
elements $(x^-_{1,r})^{k_1}\cdots
(x^-_{r,r})^{k_r}v_{\xi^\lambda}$, $k_1,\cdots ,k_r\in \bn$, we
find easily that $\eta_1=\xi_1^\lambda-k_1 -k $ for some $k\ge 0$.
Hence $k_1\le \xi^\lambda_1 - \eta_1 \le \xi^\lambda_1
-\xi^\lambda_2=m_1$ and the lemma is proved.
\end{pf}

\subl{} We can now prove the following stronger statement.
\begin{prop}\label{prop_fd}
Let $\lambda = \sum m_i \omega_i\in P^+$. Then we have,
\begin{equation} V(\lambda)=\sum_{\{k_i: 0\le k_i\le m_i, 1\le i\le r\}}
\bu(\lie n^-_{r-1})(x^-_{1,r})^{k_1}\cdots
(x^-_{r,r})^{k_r}v_\lambda.\end{equation}
\end{prop}

\begin{proof}
 Let $\lie n_{m,n}^\pm$ be the subalgebra of $\lie g$ spanned by
$x_{i,j}^\pm$, $m\le i\le j\le n$, let $\lie h_{m,n}$ be the
subalgebra spanned by  $h_i$, $m\le i \le n$ and, finally, set
$\lie g_{m,n} = \lie n_{m,n}^+ \oplus \lie h_{m,n} \oplus \lie
n_{m,n}^-$. Clearly, $\lie g_{m,n}$ is isomorphic to ${\lie
sl}_{n-m+2}$. We proceed by induction on $r$. For $r=1$, the
result follows from the defining relations of $V(\lambda)$. Note
that $\bu(\lie g_{2,r})v_\lambda\cong V(\lambda|_{\lie h_{2,r}})$
as $ \lie g_{2,r}$--modules, and hence we have by the induction
hypothesis
\begin{equation*}
\bu(\lie n^-_{2,r})v_\lambda=\sum_{\{k_i: 0\le k_i\le m_i,  2\le
i\le r, \}} \bu(\lie n^-_{2,r-1}) (x^-_{2,r})^{k_2}\cdots
(x^-_{r,r})^{k_r}v_\lambda.
\end{equation*}
Combining this  with Lemma~\ref{lem_fd}, we find that
$$V(\lambda)=\sum_{\{k_i: 0\le k_i\le m_i, 1\le i\le r\}} \bu(\lie
n^-_{r-1}) (x^-_{1,r})^{k_1} \bu(\lie n^-_{2,r-1})
(x^-_{2,r})^{k_2}\cdots (x^-_{r,r})^{k_r}v_\lambda.$$  Since $\lie
n^-_{2,r-1} \subset \lie n^-_{r-1}$ and $[x^-_{1,r},\lie
n^-_{2,r-1}]=0$, the proposition follows.
\end{proof}

\subsection{Some results on $\bu(\lie{sl}_2[t])$ and their  consequences}

 \subl{} Assume in this section that $\lie g=\lie{sl}_2$ and set
$\omega = \omega_1$, $x^\pm = x^\pm_{1,1}$. For $n\in \bn$, let
$I_n\subset\bu(\lie n^-[t])$ be the ideal generated by the
elements $\bpr((x^+ \otimes t)^k (x^-\otimes 1)^l) $, where $k\ge
0$ and $ l>n.$

\begin{thm}\label{thm_sl2} \cite[Section 6]{CPweyl}.
\begin{enumerit}
\item
The elements  $\{\bx^-(\ell,\bos)w_{n\omega}: (\ell,\bos)\in
\bF(n)\}$ are a basis of $W(n\omega)$.
\item
The  map $\bu(\lie n^-[t]) \to W(n\omega)$ sending $g$
 to  $gw_{n\omega}$ induces an isomorphism of vector spaces
  $$\bu(\lie n^-[t])/I_n\cong W(n\omega).$$.
\end{enumerit}
\hfill\qedsymbol
\end{thm}
 For $n\in\bn$, let $ J_n$ be the subspace of $\bu(\lie n^-[t])$ spanned by
the set $$\{\bpr(\bx^+(\ell,\bos) (x^- \otimes 1)^{m}: m>n, \
(\ell,\bos)\in\bF_+\}.$$  Note that  if $m>n$, we have
$(x^-\otimes
 t^m)w_n=0$ and hence
$$\bx^+(\ell,\bos) (x^- \otimes 1)^{m} w_n=0,$$ which is
equivalent  by Lemma~\ref{prl} to
 $$\bpr(\bx^+(\ell,\bos) (x^-\otimes 1)^m) w_n=0,$$ i.e $\bpr(\bx^+(\ell,\bos) (x^-\otimes
 1)^m)\in I_n$ and hence $J_n\subset I_n$. In what follows, we
 shall prove the reverse inclusion.

\subl{} We begin with,

\begin{lem}\label{lem_ideal} Let $m,s\in\bn$ and assume that
$m>n$.
\begin{enumerit}
\item  We have
$(x^- \otimes t^s) (x^-\otimes
 1)^m\in J_n $.
\item   We have
$\bpr((h \otimes t^s) (x^-\otimes
 1)^m)\in J_n $.
\item For $g \in J_n$, $(\ell, \bos) \in \bF_+$ we have $\bpr( \bx^+(\ell,\bos) g) \in J_n$.
\end{enumerit}
\end{lem}

\begin{proof}
Note that  $$[x^+ \otimes t^s, (x^-\otimes 1)^{m+2}]= -(m+2)(m+1)
(x^- \otimes t^s) (x^-\otimes 1)^{m}. $$  Proposition
\ref{proppr}~(ii) gives, $$\bpr((x^+ \otimes t^s)(x^-\otimes
1)^{m+2})=\bpr([x^+ \otimes t^s, (x^-\otimes 1)^{m+2}])=
-(m+2)(m+1) (x^- \otimes t^s) (x^-\otimes 1)^{m},$$and hence  (i)
follows.  The proof of (ii) is similar, using $$[h\otimes t^s
,(x^-\otimes 1)^{m}] = - 2m (x^- \otimes t^s) (x^-\otimes 1)^{m-1}
= \frac{2}{m+1} \bpr((x^+ \otimes t^s) (x^-\otimes 1)^{m+1}).$$
 Next, note that
Proposition~\ref{proppr}~(iii) gives, $$\bpr(
\bx^+(\ell,\bos)\bpr(\bx^+(\ell',\bos') (x^- \otimes 1)^m )) =
\bpr(\bx^+(\ell,\bos)\bx^+(\ell',\bos')(x^- \otimes 1)^m).$$ So
Part~(iii) follows since
$\bx^+(\ell,\bos)\bx^+(\ell',\bos')=\bx^+(\ell'',\bos'')$ where
$\ell''=\ell+\ell'$ and $\bos''$ is obtained by  concatenating
$\bos$ and $\bos'$ into a partition.
\end{proof}

\subl{}  Note that $J_n$ contains the generators of $I_n$ and so
if we show that $J_n$ is an ideal we establish that $J_n = I_n$.
\begin{prop}\label{prop_id} For all $n\in\bn$, the subspace
$J_n$ is an ideal in $\bu (\lie n^-[t])$ and hence $J_n=I_n$.
\end{prop}
\begin{proof} Let $k,m\in\bn$ and $m>n$. Then,
\begin{eqnarray*}&(x^-\otimes t^k)\cdot \bpr(\bx^+(\ell,\bos) (x^- \otimes
1)^{m})&=\bpr((x^-\otimes t^k)\bx^+(\ell,\bos) (x^- \otimes
1)^{m}),\\&&=\bpr(\bx^+(\ell,\bos)(x^-\otimes t^k)(x^- \otimes
1)^{m})\\ &&+\bpr([x^-\otimes t^k,\bx^+(\ell,\bos)] (x^- \otimes
1)^{m}).\end{eqnarray*} Lemma~\ref{lem_ideal} (i), (iii) implies
that
 $\bpr(\bx^+(\ell,\bos)(x^-\otimes t^k)(x^- \otimes
1)^{m}) \in J_n$. For the second term, on the right hand side of
the preceding equation, note
that $[x^-\otimes t^k,\bx^+(\ell^+,\bos^+)]$ is a linear
combination of terms of the form $\bx^+(\ell',\bos')$ and
$\bx^+(\ell',\bos')(h\otimes t^p)$ where $(\ell', \bos') \in
\bF_+$, $p>0$. Now, if  $m>n$, the element
$\bpr(\bx^+(\ell',\bos') (x^- \otimes 1)^{m}) \in J_n$ by
definition . By Lemma~\ref{lem_ideal} (ii), (iii) we have
$$\bpr\left(\bx^+(\ell',\bos') (h\otimes t^p) (x^- \otimes
1)^{m}\right) = \bpr \left(\bx^+(\ell',\bos') \bpr\left((h\otimes
t^p) (x^- \otimes 1)^{m}\right)\right)\in J_n.$$  This proves that
$(x^-\otimes t^k)\cdot \bpr(\bx^+(\ell,\bos) (x^- \otimes
1)^{m})\in J_n$ establishing the proposition.
\end{proof}

\subl{} The following is now an immediate corollary of
Theorem~\ref{thm_sl2} and  Proposition~\ref{prop_id}.

\begin{prop}\label{bas}
 Assume that $\lie g$ is of type $\lie{sl}_2$ and
that $n\ge 0$.  Then, $\bu(\lie n^-[t])$ is spanned by elements of
the set  $$\left\{\bx^-(\ell,\bos):(\ell,\bos)\in
\bF(n)\,\right\}\,\bigcup\,
 \left\{\bpr(\bx^+(\ell,\bos)(x^- \otimes 1)^{m}): m>n,\
(\ell,\bos)\in\bF_+\right\}.$$\hfill\qedsymbol
\end{prop}

 We now return to the case of $\lie{sl}_{r+1}$.  Fix $1\le i\le r$.
Applying  Proposition~\ref{bas} to the subalgebra of $\lie{sl}_{r+1}$
generated by the elements $x_{i,r}^\pm$,   and using
Proposition~\ref{prop_prg}, we obtain the following result.

\begin{cor}\label{prop_sl2} Let $1\le i\le r$ and $n\ge 0$.
For all $(\ell,\bos)\in\bF$, the element  $\bx^-_{i,r}
(\ell,\bos)$ is in the span of the union of
$$\left\{\bx^-_{i,r}(\ell,\bos'):(\ell,\bos')\in \bF(n), \
|\bos'|=|\bos|\right\},$$ and $$
 \left\{\bpr(\bx^+_{i,r}(m-\ell,\bos')(x^-_{i,r} \otimes 1)^{m}): m>n,\
(m-\ell,\bos')\in\bF_+,\ \  |\bos'|=|\bos| \right\}.$$ In
particular, the element $\bx_{i,r}^-(\ell,\bos)$, is in the span
of
 $$
 \left\{\bpr(\bx^+_{i,r}(m-\ell,\bos')(x^-_{i,r} \otimes 1)^{m}): m\ge 0,\
(m-\ell,\bos')\in\bF_+, \ \  |\bos'|=|\bos| \right\}.$$
\hfill\qedsymbol
\end{cor}

\subsection{A crucial Lemma}

\subl{}

The goal of this section is to prove the following statement

\begin{lem}\label{lem_dif}  Let $(\ul^+\us^+)\in \bF^r_+$. For all $(\ul,\us)\in
\bF^r$ with $\ell_i\ge \ell_i^+$, $1\le i\le r$, we have
\begin{equation}\label{eq_dif}
\bpr\left(\bxu^+ (\ul^+, \us^+) \bxu^- (\ul, \us)\right) -
\prod_{i=1}^r \bpr \left( \bx^+_{i,r} (\ell^+_i, \bos^+_i)
\bx^-_{i,r} (\ell_i, \bos_i)\right) \in \bu(\lie n^-[t])^{>(\ul -
\ul^+, |\us| + |\us^+|)}
\end{equation}
In particular, $$\bpr\left(\bxu^+ (\ul^+, \us^+) \bxu^- (\ul,
\us)\right) \in \bu(\lie n^-[t])^{\ge(\ul - \ul^+, |\us| +
|\us^+|)}.$$

\end{lem}

We shall  prove the first statement by induction on $|\ul^+| =
\ell^+_1 + \dots +\ell^+_r$. Proposition~\ref{prop_prg} gives
$$\prod_{i=1}^r \bpr \left( \bx^+_{i,r} (\ell^+_i, \bos^+_i)
\bx^-_{i,r} (\ell_i, \bos_i)\right) \in \bu(\lie n^-[t])^{(\ul -
\ul^+, |\us| + |\us^+|)},$$ which proves the second statement of
the proposition.

\subl{} We shall need the following. For a partition $\bos =
(\bos(1) \le \dots \le \bos(\ell))$ let $$\bos^{(i)} = (\bos(1)
\le \dots \le \bos(i-1) \le \bos(i+1) \le \dots \le \bos(\ell)).$$

\begin{lem}\label{comm} Let $(\ell,\bos)\in \bF$, $s\in\bn$
 and  $1\le i\le  j\le r$.  Then
$$[x^+_{i,r}\otimes t^s, \bx_{j,r}^-(\ell,\bos)] = \left\{
\begin{array}{ll}
\sum_{k=1}^\ell \bx_{j,r}^-(\ell-1, \bos^{(k)}) \left(
x^+_{i,j-1}\otimes t^{s+ \bos(k)}\right) & i<j,\\ \sum_{k=1}^\ell
\bx_{j,r}^-(\ell-1, \bos^{(k)}) \left(x^-_{j,i-1} \otimes t^{s+
\bos(k)}\right)\ & i>j,\\ \sum_{k=1}^\ell \bx_{j,r}^-(\ell-1,
\bos^{(k)}) \left( h_{i,r} \otimes t^{s+ \bos(k)}\right) +
\bpr((x^+_{i,r}\otimes t^s) \bx_{j,r}^-(\ell,\bos)) & i=j,
\end{array}
\right.$$ where $h_{i,r} = h_i + \dots + h_r \in \lie h$.
\end{lem}

\begin{proof}
 The commutators are calculated as usual by using the fact that
$\lie g[t]$ acts on $\bu(\lie g[t])$ as derivations (see
~\eqref{leibn}). The case $i<j$ follows from $$[x^+_{i,r}\otimes
t^s, x^-_{j,r}\otimes t^{s'}] = x^+_{i,j-1}\otimes t^{s+s'},
\qquad i<j,$$ and the observation that $x^+_{i,j-1}\otimes
t^{s+s'}$ commutes with the factors of $\bx_{j,r}^-(\ell,\bos)$.

The case $i>j$  follows from $$[x^+_{i,r}\otimes t^s,
x^-_{j,r}\otimes t^{s'}] = x^-_{j,i-1} \otimes t^{s+s'}, \qquad
i>j,$$ and the observation that $x^-_{j,i-1} \otimes t^{s+s'}$
commutes with the factors of $\bx_{j,r}^-(\ell,\bos)$.

Finally consider the case  $i=j$. We have $$[x^+_{i,r}\otimes
t^s, x^-_{i,r}\otimes t^{s'}] = h_{i,r} \otimes t^{s+s'},$$  which
gives \begin{eqnarray*}&[x^+_{i,r}\otimes t^s,
\bx_{i,r}^-(\ell,\bos)]& = \sum_{k=1}^\ell \left(
\prod_{n=1}^{k-1} x^-_{i,r}\otimes t^{\bos(n)} \right) (h_{i,r}
\otimes t^{s+\bos(k)}) \left( \prod_{n=k+1}^{\ell}
x^-_{i,r}\otimes t^{\bos(n)} \right)\\ && =\sum_{k=1}^\ell
\bx_{j,r}^-(\ell-1, \bos^{(k)}) \left( h_{i,r} \otimes t^{s+
\bos(k)}\right)\\&& + \sum_{k=1}^\ell \left( \prod_{n=1}^{k-1}
x^-_{i,r}\otimes t^{\bos(n)} \right) \left[h_{i,r} \otimes
t^{s+\bos(k)}, \left( \prod_{n=k+1}^{\ell} x^-_{i,r}\otimes
t^{\bos(n)} \right)\right].
\end{eqnarray*}
The first term on the right hand side of the preceding  equality
is in $\bu(\lie n^-[t])\lie h[t]$. Since  $[\lie h[t],\lie
n^-[t]]\subset\lie n^-[t]$ the second term is in $\bu(\lie
n^-[t])$ and hence  is equal to  $\bpr([x^+_{i,r}\otimes
t^s,\bx_{i,r}^-(\ell,\bos)])$. Since the latter term is  also
equal to  $ \bpr((x^+_{i,r}\otimes t^s)\bx_{i,r}^-(\ell,\bos))$,
the proof of the Lemma is complete. \end{proof}

\subsubsection{The case $|\ul^+|=1$ in the proof of Lemma
\ref{lem_dif}}
 In this case, we have $\bxu^+ (\ul^+, \us^+) = x_{i,r}^+ \otimes t^s$
  for some $1\le i\le r$ and $s>0$.
 For $1\le j\le r$, set
 $$g_j=\bx^-_{1,r}(\ell_1,\bos_1)\cdots\bx^-_{j-1,r}(\ell_{j-1},\bos_{j-1}),\
 \
 g_j'=\bx^-_{j+1,r}(\ell_{j+1},\bos_{j+1})\cdots\bx^-_{r,r}(\ell_{r},\bos_{r}),$$
 where we understand that $g_1=1$ and $g_r'=1$.  Using~\eqref{leibn} we get,
\begin{eqnarray*}&\bpr\left( (x_{i,r}^+ \otimes t^s) \cdot \bxu^-(\ul, \us)
\right) &= \bpr\left([x_{i,r}^+ \otimes t^s, \bxu^-(\ul,
\us)]\right)\\ && =\sum_{j=1}^r\bpr\left(g_j [x_{i,r}^+ \otimes
t^s, \bx^-_{j,r} (\ell_j, \bos_j)]g_j'\right).\end{eqnarray*}

If $i<j$, then by using Lemma \ref{comm} we see that
$$\bpr\left(g_j [x_{i,r}^+ \otimes t^s, \bx^-_{j,r} (\ell_j,
\bos_j)]g_j'\right)=g_j\sum_{k=1}^{\ell_j}\bpr\left(
\bx_{j,r}^-(\ell_j-1, \bos_j^{(k)})( x^+_{i,j-1} \otimes t^{s+
\bos_j(k)})g_j'\right)=0,$$ where the last equality follows
 since $x^+_{i,j-1} \otimes t^{s+\bos_j(k)}$ belongs to $\lie
b^+[t]$ and commutes with $g_j'$.

If $i>j$, then by Lemma \ref{comm}, we get
 \begin{eqnarray*}&\bpr\left(g_j
[x_{i,r}^+ \otimes t^s, \bx^-_{j,r} (\ell_j,
\bos_j)]g_j'\right)&=g_j\sum_{k=1}^{\ell_j} \bx_{j,r}^-(\ell_j-1,
\bos_j^{(k)})( x^-_{j,i-1} \otimes t^{s+ \bos_j(k)})g_j',\\ &&=(
x^-_{j,i-1} \otimes t^{s+ \bos_j(k)})g_j\sum_{k=1}^{\ell_j}
\bx_{j,r}^-(\ell_j-1, \bos_j^{(k)})g_j',\end{eqnarray*} where the
second equality follows since  $x^-_{j,i-1} \otimes t^{s+
\bos_j(k)}$ commutes with both $\bx_{j,r}^-(\ell_j-1,\bos_j^{(k)})$ and
$g_j$. Moreover,   the right hand side of the second equation
clearly lies in $\sum_{k=1}^{\ell_i}\bu(\lie n^-[t])^{\boi_k}$
where
$$\boi_k=(\ul-\boe_j,|\us|-\bos_j(k)\boe_j)>(\ul-\boe_i,
|\us|+s\boe_i).$$

Finally, if  $i=j$ we have
\begin{eqnarray*}
&\bpr\left(g_i [x_{i,r}^+ \otimes t^s, \bx^-_{i,r} (\ell_i,
\bos_i)]g_i'\right)&=g_i \sum_{k=1}^{\ell_i} \bx^-_{i,r}(\ell_i
-1, \bos_i^{(k)}) \bpr\left((h_{i,r} \otimes t^{s+ \bos_i(k)})
g_i'\right) \\ && + g_i\bpr\left(\bpr\left((x^+_{i,r}\otimes t^s)
\bx_{i,r}^-(\ell,\bos)\right)g_i'\right)
\end{eqnarray*}
It is now not hard to see that the first term is in $$
\sum_{i<m\le r} \bu(\lie n^-[t])^{> (\ul - \boe_i, |\us| -
\bos_i(k) \boe_i+ (s+\bos_i(k)) \boe_m)} \subset  \bu(\lie
n^-[t])^{> (\ul - \boe_i, |\us| + s\boe_i)}.$$   And the second
term is equal to $g_i \bpr\left((x^+_{i,r}\otimes t^s)
\bx_{i,r}^-(\ell,\bos)\right)g_i'$ since  both
$\bpr\left((x^+_{i,r}\otimes t^s) \bx_{i,r}^-(\ell,\bos)\right)$
and $g_i'$ belongs to $\bu (\lie n^-[t])$.

\subsubsection{The inductive step}

For the inductive step,  fix $1\le i_0\le r$  such that
$\ell^+_{i_0}\ne 0$. Then  we can write, $\bxu^+ (\ul^+, \us^+) =
(x^+_{i_0,r}\otimes t^s)\bxu^+ (\ul', \us')$ for $\ul' = \ul -
\boe_{i_0}$ and some $s\in\bn$, $\us'\in \bF^r$. Note that by
Proposition~\ref{proppr}~(iii) we have,
$$\bpr\left((x^+_{i_0,r}\otimes t^s)\bxu^+ (\ul', \us')
\bxu^-(\ul, \us)\right) = \bpr\left((x^+_{i_0,r}\otimes
t^s)\bpr\left(\bxu^+ (\ul', \us') \bxu^-(\ul, \us)\right)
\right).$$ The result follows by first using the induction
hypothesis for  $\bpr\left(\bxu^+ (\ul', \us')
\bxu^-(\ul,\us)\right)$  and then the result when $|\ul^+| = 1$.

\subsection{Proof of Proposition~\ref{assgraded}~(ii)}.

 \subl{} To complete the proof of the proposition, we need some
 more results. We begin with  two lemmas  related to the order.

\begin{lem}\label{lem_l} Let $\ul,\, \ul' \in \bn^r$.
\begin{enumerit}
\item We have $\eta_r(\ul)=\eta_r(\ul')$ if and only if $\ul=\ul'$.
\item If $\eta_r(\ul) - \eta_r(\ul') \in Q_+$, but $\eta_r(\ul) \ne  \eta_r(\ul')$
then $\ell' \eql \ell$.
\end{enumerit}
\end{lem}

\begin{proof}
Part (i) follows since $\alpha_{i,r}$, $1\le i\le r$ are linearly
independent. For part (ii),  write  $\eta_r(\ul) - \eta_r(\ul') =
\sum_{i=1}^r m_i \alpha_{i}$, where
$m_i=\sum_{j<i}(\ell_j-\ell_j')$. Since $m_i\ge 0$ for all $1\le
i\le r$,  we can  take $1\le s\le r$ minimal such that $m_s>0$.
This means that   $\ell'_i = \ell_i$ for $i<s$ and $\ell'_s <
\ell_s$, so the lemma is proved.
\end{proof}

\subl{}

 \begin{lem}\label{lem_d}  Assume that
  $g \in \bu(\lie n^-_{r-1}[t])_{-\eta}$, $\eta \ne 0$
 and let $(\ul, \us) \in \bF^r_+$.
Then, $[\bxu^+ (\ul, \us), g]$ is a linear combination of elements
of the form $g' \bxu^+ (\ul', \us')$ where  $g' \in \bu(\lie
n^-_{r-1}[t])$, $(\ul', \us') \in \bF^r_+$ and $\us' \eqd \us.$
\end{lem}

\begin{proof}
It suffices to prove the result when $g$ is a product of terms of
the form $x_{j,k}^-\otimes t^s$, where $1\le j\le k\le r-1$ and
$s\ge 0$. We prove this by induction on the number $m$ of terms in
the product.
 If $g = x^-_{j,k} \otimes t^s$, $1\le j\le k<r$, then
$[x^+_{i,r} \otimes t^{s'}, g] = - \delta_{i,j} x^+_{k+1,r}
\otimes t^{s+s'}$. It follows from~\eqref{leibn} that $[\bxu^+
(\ul, \us), g]$ is a linear combination of terms of the form $
\bxu^+ (\ul', \us')$ where $(\ul', \us') \in \bF^r_+$ and $$|\us'|
= |\us| - \bos_j(m)\boe_j + (\bos_j(m)+s)\boe_{k+1},\quad 1\le
m\le \ell_j.$$ Since $\bos_j(m)>0$ for $1\le m\le \ell_j$, we have
$\us' \eqd \us$, hence induction begins.

For the inductive step, write   $g = (x^-_{j,k} \otimes t^s)g_0$,
where $g_0$ is a product of $(m-1)$--elements of the form
$x_{j,k}^-\otimes t^{s'}$, where $1\le j\le k\le r-1$ and $s'\ge
0$.  We have, $$[\bxu^+ (\ul, \us), (x^-_{j,k} \otimes t^s)g_0] =
[\bxu^+ (\ul, \us), (x^-_{j,k} \otimes t^s)]g_0 + (x^-_{j,k}
\otimes t^s)[\bxu^+ (\ul, \us), g_0].$$ The second term has the
required form by the induction hypothesis. For the first one
 we have (by using the induction hypothesis for $g=x_{j,k}^-\otimes t^s$) that
$$[\bxu^+ (\ul, \us), x^-_{j,k} \otimes t^s]g_0 =
\sum_{{\ul',\us'}\atop{\us'\eqd\us}} a_{\ul',\us'}\bxu^+(\ul',
\us') g_0 =
 \sum_{{\ul',\us'}\atop{\us'\eqd\us}} a_{\ul',\us'}\ g_0\bxu^+(\ul', \us') +
\sum_{{\ul',\us'}\atop{\us'\eqd\us}} a_{\ul',\us'}[\bxu^+(\ul',
\us'), g_0],$$ where $a_{\ul',\us'} \in \bc$.
 The result follows by using the induction hypothesis again.
\end{proof}

\subl{}
\begin{prop}\label{prop_bound1} Let
$\lambda=\sum_{i=1}^rm_i\omega_i\in P^+$.  Fix $1\le k\le r$ and
$n\in\bn$ with $n>m_k$. Let $(\ell_i, \bos_i) \in \bF$, $1\le i
\ne k \le r$, and $(n-\ell_k, \bos_k) \in \bF_+$.  Then,
\begin{equation*}
\prod_{i=1}^{k-1} \bx^-_{i,r}(\ell_i, \bos_i)
\bpr\left(\bx_{k,r}^+(n-\ell_k,\bos_k) (x_{k,r}^- \otimes
1)^{n}\right)  \prod_{i=k+1}^{r} \bx^-_{i,r}(\ell_i, \bos_i)
w_\lambda \in  W(\lambda)^{> (\ul, |\us|)}
\end{equation*}
\end{prop}

\begin{proof}
Applying Corollary~\ref{prop_sl2}  for all $1\le i \le r$ except
$i=k$, we see that it suffices to prove that  elements of the form
$$ \prod_{i=1}^r\bpr\left(\bx_{i,r}^+(n_i-\ell_i,\bos^+_i)
(x_{i,r}^- \otimes 1)^{n_i}\right)w_\lambda, \qquad
(n_i-\ell_i,\bos_i^+)\in{\bF}_+,\ \   n_i\ge \ell_i,\ \ 1\le i\le
r,$$
 are in $W(\lambda)^{>(\ul, |\us^+|)}$ if $n_k > m_k$.
By Lemma~\ref{lem_dif} and  Lemma~\ref{prl} this is equivalent to
proving that $$w_1=\prod_{i=1}^r \bx^+_ {i,r} (n_i-\ell_i,
\bos^+_i) \prod_{i=1}^r (x^-_{i,r}\otimes 1)^{n_i} w_\lambda\in
W(\lambda)^{>(\ul,\, |\us^+|)}.$$ Note that $w_1\in
W_{\eta_r(\ul)}$. Now, Proposition~\ref{prop_fd} implies that
$w_1$ is a linear combination of elements from $W_{\eta_r(\ul)}$
of the form $$w_2=\prod_{i=1}^r \bx^+_{i,r} (n_i-\ell_i, \bos^+_i)
g \prod_{i=1}^r (x^-_{i,r}\otimes 1)^{n_i'} w_\lambda$$ with
$n_k'\le m_k<n_k$ and $g \in \bu(\lie n^-_{ r-1} \otimes 1)$. This
gives $g \in \bu(\lie n^-)_{-\nu}$, where   $$\nu = \sum_{j=1}^{r}
(n_j - n_j') \alpha_{j,r}. $$  Since $n_k'<n_k$ and
$\alpha_{i,r}$, $i=1\dots r$, are linearly independent, we see
that $\nu \ne 0$, in particular, $g$ is not a constant. Now, Lemma
\ref{lem_d} implies that $w_2$ is a linear combination of elements
from $W_{\eta_r(\ul)}$ of the form
  $$w_3=g' \prod_{i=1}^r
\bx^+_{i,r} (n_i - \ell_i', \bos_i') \cdot \prod_{i=1}^r
(x^-_{i,r}\otimes 1)^{n_i'} w_\lambda, \qquad g' \in \bu(\lie
n^-_{r-1}[t]),$$ where either $g' = g$, $\ul' = \ul$, $\us' =
\us^+$
 or $|\us'| \eqd |\us^+|$.

By Lemma~\ref{lem_dif} we have $$ \prod_{i=1}^r \bx^+_{i,r} (n_i -
\ell_i', \bos_i') \cdot \prod_{i=1}^r (x^-_{i,r}\otimes 1)^{n_i'}
\in \bu(\lie n^-[t])^{\ge (\ul'',\,|\us'|)}, \quad \mbox{where} \
\ \ell''_i = \ell_i' - n_i + n_i',$$ therefore $w_3 \in
W(\lambda)^{\ge (\ul'',\,|\us'|)}$, and hence the proposition
follows if we show that $(\ul'',\,|\us'|) > (\ul, |\us^+|)$. For
this, observe that since  $w_3 \in W(\lambda)_{\eta_r(\ul)}$, we
have $g' \in \bu(\lie n^-[t])_{-\nu'}$, where $\nu'= \eta_r(\ul) -
\eta_r(\ul'')$. Since $\nu' \in Q^+$, by Lemma~\ref{lem_l} we have
either $\ul'' \eql \ul$ or $\nu' =0$, $\ul'' = \ul$. In the first
case which includes the case $g' = g$, we have  $(\ul'',\,|\us'|)
> (\ul, |\us^+|)$ by the definition of the order, and in the the
second case we have $\ul'' = \ul$, but $|\us'| \eqd |\us^+|$, so
$(\ul'',\,|\us'|) > (\ul, |\us^+|)$ and the proof is complete.
\end{proof}

\subl{} {\it Proof of Proposition~\ref{assgraded}~(ii)}.  Let
$\lambda=\sum_{i=1}^r m_i\omega_i\in P^+$. Using
Corollary~\ref{prop_sl2}  simultaneously for $1\le i\le r $ and
for $n = m_1, \dots, m_r$, we see that $\bu(\lie
n^-[t])^{(\ul,\ud)}$ is spanned by elements from
\begin{equation}\label{these}
g\bxu^-(\ul,\us), \qquad g\in\bu(\lie n_{r-1}[t]),\ \ |\us| = \ud,
\ \
(\ul,\us) \in \bF^r(\lambda)
\end{equation}
 together with the elements from
\begin{equation}\label{other}
g \prod_{i=1}^{k-1}\bx^-_{i,r}(\ell_i,\bos_i)\cdot
\bpr((\bx^+_{k,r}(m-\ell_k,\bos)(x_{k,r}^-\otimes 1)^m)
\prod_{i=k+1}^{r}\bx^-_{i,r}(\ell_i,\bos_i), \qquad m> m_k, \ \,
\ \ |\us| = \ud,
\end{equation}
where $g\in\bu(\lie n_{r-1}[t])$,  $(\ell_i, \bos_i) \in \bF$,
$1\le i \ne k \le r$, and $(m-\ell_k, \bos_k) \in \bF_+$.
Proposition~\ref{prop_bound1} implies that the elements
of~\eqref{other} applied to to $w_\lambda$ lie in $W(\lambda)^{>
(\ul,\ud)}$ and hence, $W(\lambda)^{\ge (\ul,\ud)}/ W(\lambda)^{>
(\ul,\ud)}$ is spanned by the elements in~\eqref{these} applied to
$w_\lambda$. In turn, this means that $\gr W(\lambda)$ is spanned
by $g\bxu^-(\ul,\us) w_\lambda$ with $(\ul,\us) \in
\bF^r(\lambda)$, $g\in\bu(\lie n_{r-1}[t])$,
that is by the images of $\psi^{\ul,\us}$ for such $(\ul,\us)$ and
the proof is complete.

\section*{ Index of Notation} We provide for the readers convenience a brief index of the notation
which is used repeatedly in this paper.

 \vskip 12pt

Section 1.1.1: $\bz$, $\bn$, $\bn_+$,
$\lie g$, $\lie h$, $\lie n^\pm$, $x_{i,j}^\pm$, $H_i$, $h_i$,
$\omega_i$, $\alpha_i$, $\alpha_{i,j}$, $P$, $P^+$, $Q$, $Q^+$, $\bu(\lie a)$.

\vskip 12pt

 Section 1.1.2:  $\bz[P]$, $e(\mu)$, $\ch_{\lie g}(V)$, $V(\lambda)$, $v_\lambda$.

\vskip 12pt

 Section 1.1.3:  $\lie a[t]$, $\lie at[t]$, $\ch_t(M)$.

\vskip 12pt

 Section 1.2.1:  $W(\lambda)$, $w_\lambda$.

\vskip 12pt

 Section 1.2.2: $W(\lambda)_s$, $W(\lambda)_\mu$, $W(\lambda)_{\mu,s}$.

\vskip 12pt

 Section 2.1.1:  $\bF$, $(\ell, \bos)$, $(0, \emptyset)$,
$\bx_{i,j}^\pm(\ell,\bos)$, $(\ul,\us)$, $\underline{\bx}_j^\pm(\ul,\us)$, $\cal B^r$.

 \vskip 12pt

 Section 2.1.2:
$\bF(m)$, $\bF^j(\lambda)$, $\eta_j(\ul)$, $\cal B^r(\lambda)$.

\vskip 12pt

 Section 2.2.1: $\lie g_{r-1}$, $\lie n^\pm_{r-1}$,
$\lie u^\pm_r$.

\vskip 12pt

Section 2.2.2: $\bi^r$, $\eqd$, $\eql$, $>$, $|\bos|$, $|\us|$,
$W(\lambda)^{\ge\boi}$, $W(\lambda)^{>\boi}$,  $\gr( W(\lambda))$,
$\gr^\boi$.

\vskip 12pt

Section 2.2.3:  $\psi^{\ul,\us}$.

\vskip 12pt

Section 3.1.1: $\lie b^+$, $\bpr$.

\vskip 12pt

Section 3.1.2: $\bu(\lie n^-[t])^\boi$,  $\bu(\lie n^-[t])^{\ge
\boi}$,  $\bu(\lie n^-[t])^{>\boi}$.

\vskip 12pt

Section 3.1.3: $\bu(\lie g[t])_{\mu}$, $\bu(\lie g[t])_{\mu,s}$,
 $\bu(\lie n^\pm[t])_{\mu}$, $\bu(\lie n^\pm[t])_{\mu,s}$.

\vskip 12pt

Section 3.1.4:  $\bF_+$.

\vskip 12pt

Section 3.2.1:  $\boe_i$.

\end{document}